\documentclass[11pt]{article}

\usepackage{amsthm,amsfonts,amssymb,amsmath,color}

\numberwithin{equation}{section}

\renewcommand\b{\beta}
\renewcommand\o{\omega}

\def\g{\gamma}

\def\O{\Omega}
\def\th{\theta}

\def\l{\lambda}

\def\epsilon{\varepsilon}
\def\e{\varepsilon}

\newcommand\br{\begin{rem}}
\newcommand\er{\end{rem}}
\newcommand\bp{\begin{pmatrix}}
\newcommand\ep{\end{pmatrix}}
\newcommand\be{\begin{equation}}
\newcommand\ee{\end{equation}}
\newcommand\ba{\begin{equation}\begin{aligned}}
\newcommand\ea{\end{aligned}\end{equation}}

\newcommand\nn{\nonumber}

\newcommand{\ZZ}{{\mathbb Z}}

\newcommand{\Id}{{\rm Id }}

\newcommand{\tr}{{\rm tr }}
\newcommand{\supp}{{\rm supp }}

\newcommand{\diag}{\text{\rm diag}}

\newtheorem{defi}{Definition}[section]
\newtheorem{theorem}[defi]{Theorem}

\newtheorem{lemma}[defi]{Lemma}

\newtheorem{remark}[defi]{Remark}

\newtheorem{assumption}[defi]{Assumption}

\def\op{{\rm op} }

\numberwithin{equation}{section}
\hoffset - 1 cm

\begin{document}

\title{A note on separation conditions of resonance sets in the instability analysis for high-frequency oscillations in geometric optics}

\author{ Jiaojiao Pan\footnote{Corresponding author, Department of Mathematics, Nanjing University, 22 Hankou Road, Gulou District, Nanjing 210093, China, {\tt panjiaojiao@smail.nju.edu.cn.}}}

\date{}

\maketitle

\begin{abstract}
In this paper, we study the instability of highly-oscillating solutions to semi-linear hyperbolic systems. A instability criterion was given in \cite{Lu} under rather strong separation conditions of resonance sets: coupled resonance sets are pairwise disjoint. Here we show that such separation conditions in \cite{Lu} can be relaxed: one of the coupled non-transparent resonance sets is allowed to intersect with at most two others.  We obtain the same instability criterion as in \cite{Lu}.  Finally, we give some applications to coupled Klein-Gordon systems with equal masses and nonlinear terms specified particularly,  where on the intersections of resonance sets, the related interaction coefficients are non-transparent.

\end{abstract}




\renewcommand{\refname}{References}
\section{Introduction}
\label{Introduction}
\subsection{Problem formulation}
\label{Problem formulation}
We study the instability of highly-oscillating solutions to semilinear  systems of the form
\be\label{(1.1)}
\partial_{t}u+\frac{1}{\e}A_{0}u+\sum_{1\leq j\leq d}A_{j}\partial_{x_{j}}u=\frac{1}{\sqrt{\e}}B(u,u),
\ee
in the small wavelength limit $\e\rightarrow 0$. Here the unknown $u=u(t,x)=(u_{1},u_{2},\cdots, u_{N})\in \mathbb R^{N}$, with time variable $t\in \mathbb R_{+}$, space variable $x\in\mathbb R^{d}$. Furthermore, $A_{0}\in \mathbb R^{N\times N}$ is skew-symmetric,  $A_{j}\in \mathbb R^{N\times N}$ is symmetric, the source term $B(u,u)$ is symmetric bilinear, whose specific form can be derived from a phenomenological description of nonlinear interactions \cite{Nishiura}.

The problem mainly tells the propagation of light with relevant initial highly oscillating condition of the form
\be\label{(1.2)}
u(\e,0,x)=a(x)e^{ik\cdot x/{\e}}+ a(x)^{*}e^{-ik\cdot x/{\e}}+ O (\sqrt{\e}),
\ee

where $a$ is of high Sobolev regularity, $k$ is a given wave-number in ${\mathbb{R}}^{d}$, and $a^{*}$ denotes the complex conjugation of  $a$.   For fixed $\e>0$, the classical theory on symmetric hyperbolic systems gives the existence and uniqueness of local-in-time solutions to \eqref{(1.1)}-\eqref{(1.2)} in Sobolev spaces $H^{s}$ with $s>d/2$ and the a priori existence time is $O(\e^{1/2})$. 

Under certain compatibility condition, an approximate solution for \eqref{(1.1)}-\eqref{(1.2)} over time interval $[0,T]$ with $T>0$ independent of $\e$ can be constructed via WKB expansion. A main concern is the stability of such WKB solutions. We say a WKB approximate solution $u_{a}$ is stable if the real solution $u$ stays close to $u_{a}$ in some time interval $[0,T]$ with $T$ independent of ${\e}$, and a WKB solution $u_{a}$ is unstable provided in some short interval $[0,T_{\e}]$ with $T_{\e}\rightarrow 0$, $|u-u_{a}|$ is much amplified compared to the initial differential.

Transparency is introduced in by Joly, M\'etivier, Rauch (see \cite{Joly}, see also Dumas's survey \cite{Dumas}) which ensures the stability of WKB solution and is analogous to the null conditions which imply global existence for nonlinear wave equations,  see Klainerman's classical work \cite{Klainerman}. The link between transparency and null forms can be seen in Lannes' Bourbaki review \cite{Lannes}. 

 Joly, M\'etiver and Ranch \cite{Joly}  considered Maxwell-Bloch systems in the critical regime of geometric optics and the transparency conditions are verified. Consequently, the existence and the stability of WKB solutions to Maxwell-Bloch equations in a supercritical regime are obtained. Following Joly, M\'etivier and Rauch, it was verified by Texier  \cite{Texier2,Texier4} that the Euler-Maxwell equations satisfy a form of transparency, and by Lu that the Maxwell-Landau-Lifschitz equations also are transparent in one dimensional setting \cite{Lu1}. Cheverry, Gu\`es and M\'etivier \cite{Cheverry} showed that for systems of conservation laws, linear degeneracy of a field implies transparency. Jeanne \cite{Jeanne} showed that the Yang-Mills equations provide another example of a physical system exhibiting transparency properties.

A further question is that: if the transparency conditions are not satisfied, will the WKB approximate solutions become unstable or not? In \cite{Lu},  Texier and Lu gave a rather complete description on the stability of WKB solutions to a class of semi-linear hyperbolic systems issued from highly-oscillating initial data with large  amplitudes. They found that the stability of WKB solution is determined by the sign of some stability character index $\Gamma$ determined by the linear operator of the equation and nonlinear form of the equation. Roughly speaking, if $\Gamma>0$ which ensures the symbol of the linearized operator around a WKB approximate solution admits an eigenvalue with positive real part, the solution of the linearized system has certain exponential growth and thus is unstable; if $\Gamma<0$ which ensures the eigenvalues of the symbol of the linearized operator are pure imaginary, the solution of the linearized system stay bounded and is stable.  While, if $\Gamma = 0$ and the transparency conditions are not satisfied, the stability  analysis of the WKB solutions is more delicate, and it seems both stability and instability are possible. For example, in \cite{Lu3}, Lu considered a system of two coupled Klein-Gordon equations with different velocities and different masses and instability of the WKB solutions is discovered. With a proper choice of nonlinear source terms, it was shown in \cite{Lu3} that even the equations linearized around the leading WKB terms are initially stable (transparency conditions are satisfied), while the resonances associated with the higher-order harmonics of the WKB solutions generated by the nonlinearities can also generate instantaneous instabilities. In particular, such higher-order harmonics are not present in the data. These studies are rated to the work \cite{Lerner} and \cite{Lerner1} where the instability phenomena with loss of hyperbolicity are studied.  Later in \cite{Lu2}, Zhang and Lu considered the  non relativistic limit of the Klein-Gordon equations with quadratic nonlinearities. It is introduced compatible conditions weaker than the strong transparency conditions and a singular localization method to prove the stability of WKB solutions over long time intervals in \cite{Lu2}. In particular, such weaker compatible conditions ensure $\Gamma = 0$.

In \cite{Lu},  a key assumption is the resonance sets are essentially pairwise disjoint: related interaction coefficients are transparent in corresponding intersections of the resonance sets (see Assumption \ref{Separation conditions}). However, such reparation conditions of resonance sets may not hold in practice. For example in Section 5.2 in \cite{Lu} concerning a coupled Klein-Gordon system with equal masses and different velocities, the separation conditions are satisfied only in one dimensional setting. Actually,  the resonance sets often intersect in physical models.  Our goal in this paper is to offer certain weaker separation conditions which still ensure the main instability results in \cite{Lu} hold. This will make the abstract instability criterion theory in \cite{Lu} applicable to more physical models. 


This paper is organized as follows. In Section \ref{Problem formulation} to Section \ref{Existence of the symbolic flow $2*2$}, we mainly introduce the formulation of the problem, the existing instability criterion (Theorem \ref{A instability criterion}) and estimates of symbolic flows related to $2\times2$ block matrices for non-transparent resonance in \cite{Lu}. We give new relaxed separation conditions related to 2-coupled non-transparent resonances and 3-non-transparent resonances coupled in pairs and present new instability theorem in Section \ref{Main results}. In Section \ref{Ideas of the proof}, main ideas of the proof for Theorem \ref{New instability criterion} are presented. What's more, we show the the estimates of symbolic flows related to $3\times3$ block matrices for 2-coupled non-transparent resonances and 3-non-transparent resonances coupled in pairs respectively in Section \ref{Cases of two coupled resonances} and Section \ref{Three coupled resonances}, and obtain new instability criterion in Theorem \ref{New instability criterion}. Then we give an application in Section \ref{Coupled Klein-Gordon systems with equal masses} for new relaxed instability criterion which contains coupled Klein-Gordon systems with equal masses whose nonlinear terms are specified particularly. Finally, we show a more complete result about the instability for coupled Klein-Gordon systems with equal masses given in \cite{Lu}.
\subsection{Assumptions}\label{Assumptions}
Here we recall some basic assumptions introduced in \cite{Lu}, where the stability character is given.
\begin{assumption}\label{Smooth spectral decomposition}
 Assume that the matrix $A_{0}$ is real skew-symmetric, the matrices $A_{j}(1\leq j \leq d)$ are real symmetric, and the hermitian matrices $\{A_{0}/i+\sum_{1\leq j\leq d}\xi_{j}A_{j}\}_{\xi \in {\mathbb{R}}^{d}}$ have the spectral decomposition 
 \be\label{(1.3)}
 A_{0}/i+\sum_{1\leq j\leq d}\xi_{j}A_{j}=\sum_{1\leq j\leq J}\l_{j}(\xi)\Pi_{j}(\xi),
 \ee
where  $\l_{j}(\xi)$ are smooth eigenvalues and $\Pi_{j}(\xi)$  are smooth eigenprojectors satisfying the following bounds
\be\label{(1.4)}
|\partial_{\xi}^{\beta}\l_{j}(\xi)|\leq C_{\beta}(1+|\xi|^{2})^{(1-|\beta|)/2},\quad|\partial_{\xi}^{\beta}\Pi_{j}(\xi)|\leq C_{\beta}(1+|\xi|^{2})^{(-|\beta|)/2}, \quad\forall \beta\in \mathbb N^{d},
\ee
i.e. $\l_{j}(\cdot)\in S^{1}$ and $\Pi_{j}(\cdot)\in S^{0}.$ Here $S^{m}$,   $m\in\mathbb R$ denotes the set of matrix-valued symbols $a\in C^{\bar s}(\mathbb R_{x}^{d};C^{\infty} (\mathbb R_{\xi}^{d}))$ satisfying that for $\forall  \ \alpha,\beta\in \mathbb N^{d}$ with $|\alpha|\leq \bar s$ and some $C_{\alpha\beta}>0$, for $\forall \ (x,\xi)$,
\be\label{(1.5)}
|\partial_{x}^{\alpha}\partial_{\xi}^{\beta}a(x,\xi)|\leq C_{\alpha\beta}\langle\xi\rangle^{m-|\beta|},\quad\langle\xi\rangle:=(1+|\xi|^{2})^{\frac{1}{2}}.
\ee
\end{assumption}
\begin{assumption}\label{WKB solution}
For some $K_{a}\in \mathbb{N}$, some $T_{a}> 0$, there exists an approximate solution $u_{a}$ to \eqref{(1.1)} in $[0,T_{a}]$ satisfying
\be
\partial_{t}u_{a}+\frac{1}{\e}A_{0}u_{a}+\sum_{1\leq j\leq d}A_{j}\partial_{x_{j}}u_{a}=\frac{1}{\sqrt{\e}}B(u_{a},u_{a})+\e^{K_{a}}r_{a}^{\e}.\nn
\ee
The approximate solution has the form of a WKB expansion
\be\label{(1.6)}
u_{a}(\e,t,x)=e^{-i(k\cdot x-wt)/\e}u_{0,-1}(t,x)+e^{i(k\cdot x-wt)/\e}u_{0,1}(t,x)+\sqrt{\e}v_{a}(\e,t,x)\in {\mathbb{R}}^{N},
\ee
where the phase $(w,k)$ is a characteristic for the hyperbolic operator satisfying:
\be\label{(1.7)}
(\partial_{t}+\frac{1}{\e}A_{0}+\sum_{1\leq j\leq d}A_{j}\partial_{x_{j}} )(e^{\pm i(k\cdot x-wt)/\e}\vec{e}_{\pm})=0, \quad\vec{e}_{-1}=(\vec{e}_{1})^{*}
\ee
with $ \vec{e}_{1}$ and  $\vec{e}_{-1}$ are fixed unit vectors in ${\mathbb{C}}^{N}$, and  $u_{0,\pm1}$ are leading amplitudes polarized along $\vec{e}_{\pm1}$ i.e.
\be
u_{0,1}(t,x)=g(t,x)\vec{e}_{1},\quad u_{0,-1}(t,x)=g(t,x)^{*}\vec{e}_{-1}, \quad g\in C^{1}([0,T_{a}], H^{s_{a}}({\mathbb{R}}^{d})).\nn
\ee
At the same time, there holds $ v_{a}, r_{a}^{\e}\in C^{0}([0,T_{a}], H^{s_{a}}({\mathbb{R}}^{d}))$ and 
\be
\sup_{\e>0}(\sup_{|\alpha|\leq s_{a}}|(\e\partial_{x})^{\alpha}(v_{a},r_{a}^{\e})|_{(L^{\infty}[0,T],L^{2})}+|(\mathcal{F}(v_{a},r_{a}^{\e})|_{(L^{\infty}[0,T_{a}],L^{1})})<\infty.\nn
\ee
\end{assumption}

With certain compatibility conditions (weak transparency conditions), Assumption \ref{WKB solution}  holds true, see \cite{Joly} or \cite{Lu}. In particular,  $v_{a},r_{a}^{\e}$ are trigonometric polynomials in $\th$  with $\th:=({k\cdot x-wt})/{\e}$. Given the perturbation unknown $\dot u$ as
\be
u=:u_{a}+\dot u,\nn
\ee
then we have the perturbed system 
\be\label{(1.8)}
\partial_{t}\dot u+\frac{1}{\e}A_{0}\dot u+\sum_{1\leq j\leq d}A_{j}\partial_{x_{j}}\dot u=\frac{1}{\sqrt{\e}}B(u_{a})\dot u+\frac{1}{\sqrt{\e}}B(\dot u,\dot u)-\e^{K_{a}}r_{a}^{\e}.
\ee
where $B(\vec{u})v:=B(\vec{u},v)+B(v,\vec{u})$ is a bilinear
\begin{defi}\label{Resonances and interaction coefficients}
Given $i,j \in \{1,\cdots, J\}$, with J defined in Assumption \ref{Smooth spectral decomposition}, the ${(i,j)}$-resonance set is defined as follows
\be\label{(1.9)}
{\mathcal{R}}_{ij}:=\{\xi\in{\mathbb{R}}^{d}, \ \omega=\l_{i}(\xi+k)-\l_{j}{(\xi)}\}.
\ee
For $\xi\in{\mathbb{R}}^{d}$, the matrices
\be
\Pi_{i}(\xi+k)B(\vec{e}_{1})\Pi_{j}(\xi)\in{\mathbb{C}}^{N\times N}, \quad\Pi_{j}(\xi)B(\vec{e}_{1})\Pi_{i}(\xi+k)\in{\mathbb{C}}^{N\times N}\nn
\ee
are called $(i,j)$-interaction coefficients, and $(\l_{i}{(\xi+k)}-\l_{j}{(\xi)}-\omega)$ $(1\leq i,j \leq J)$ are called $(i,j)$-resonant phase. Furthermore, interaction coefficients $\Pi_{i}(\xi+k)B(\vec{e}_{1})\Pi_{j}(\xi)$ or $\Pi_{j}(\xi)B(\vec{e}_{-1})\Pi_{i}(\xi+k)$ is said to be transparent if for some $C>0$, there holds for all $\xi \in{\mathbb{R}}^{d}$, 
\be
|\Pi_{i}(\xi+k)B(\vec{e}_{1})\Pi_{j}(\xi)|\leq C|\l_{i}{(\xi+k)}-\l_{j}{(\xi)}-\omega|,\nn
\ee
or
\be
|\Pi_{j}(\xi)B(\vec{e}_{-1})\Pi_{i}(\xi+k)|\leq C|\l_{i}{(\xi+k)}-\l_{j}{(\xi)}-\omega|.\nn
\ee
If both $(i,j)$-interaction coefficients are transparent, then $(i,j)$-resonance is said to be transparent. Denote $\Re =\{(i,j), {\mathcal{R}}_{ij}\neq\varnothing \}$ the set containing all the non-empty resonance indices and $\Re_{0}\subset \Re$ the set of indices of which at least one of the interaction coefficients are non-transparent, i.e. $\Re_{0}=\{(i,j), {\mathcal{R}}_{ij}\neq\varnothing \ \mbox{and $(i,j)$-resonance is  non-transparent}\}.$
\end{defi}

Denote $\Gamma_{ij}$ the trace of the product of the $(i,j)$-interaction coefficients:
\be\label{(1.10)}
\Gamma_{ij}(\xi):=\tr \, \Pi_{i}(\xi+k)B(\vec{e}_{1})\Pi_{j}(\xi)B(\vec{e}_{-1})\Pi_{i}(\xi+k).
\ee
The stability index is given as follows 
\be\label{(1.11)}
{\Gamma} =
\left\{ \begin{aligned} &-1, \quad && \mbox{if} \  \max_{(i,j)\in \mathcal{R}_{0}} \sup_{\xi \in  \mathcal{R}_{ij} }\mbox{Re}\, \Gamma_{ij}(\xi)<0 \ \mbox{and} \  \max_{(i,j)\in{\mathcal{R}}_{0}} \sup_{\xi \in  \mathcal{R}_{ij} } |\mbox{Im}\,\Gamma_{ij}(\xi)|=0,\\
& 1, \quad &&\mbox{if} \ \max_{(i,j)\in{\mathcal{R}}_{0}} \sup_{\xi \in  \mathcal{R}_{ij} }\mbox{Re}\, \Gamma_{ij}(\xi)>0  \ \mbox{or} \ \max_{(i,j)\in{\mathcal{R}}_{0}} \sup_{\xi \in  \mathcal{R}_{ij} } |\mbox{Im}\, \Gamma_{ij}(\xi)|\neq 0. 
\end{aligned} \right.
\ee
The main results in \cite{Lu} say that if ${\Gamma}=1$, then with proper choice of initial data, the WKB solutions are unstable in a short time interval of order $O(\sqrt{\e} |\ln \e|)$, and  if ${\Gamma}=-1$, the  WKB solutions are stable in  time interval of order $O(1)$. 

And furthermore, we give several necessary notations
\be
\gamma_{ij}:=|\max_{\xi\in{\mathcal{R}}_{ij}}\mbox{Re}(\Gamma_{ij}(\xi)^{\frac{1}{2}})|,\nn
\ee
\be\label{(1.12)}
\gamma_{ij}^{+}:=|a|_{L^{\infty}}|\max_{\xi\in{\mathcal{R}}_{ij}^{h}}\mbox{Re}(\Gamma_{ij}(\xi)^{\frac{1}{2}})|,
\ee
\be\label{(1.13)}
\gamma^{+}:=\max_{(i,j)\in \Re_{0}}\gamma_{ij}^{+}.
\ee
where  $\mathcal R_{ij}^{h}$ is a neighborhood $\mathcal R_{ij}$ given as 
\be\label{(1.14)}
\mathcal R_{ij}^{h}:=\{\xi\in\mathbb R^{d}, |\l_{1}(\xi+k)-\l_{1}(\xi)-\omega|\leq h\}.
\ee
Given a frequency cut-off function $\chi_{ij}(\xi)\in C_{c}^{\infty}({\mathbb{R}}^{d})$ satisfying $0\leq \chi_{ij}(\xi)\leq1$, $\chi_{ij}(\xi)\equiv 1$ on the neighborhood $\mathcal R_{ij}^{h}$ of resonant set $\mathcal R_{ij}$, we use $\chi_{ij}^{\#}(\xi)$ to denote an extension of  $\chi_{ij}(\xi)$ in the sense that $(1- \chi_{ij}^{\#}(\xi)) \chi_{ij}(\xi)=0$. In a similar way, we can define  $\varphi_{ij}^{\#}(x)$ an extension of $\varphi_{ij}(x)\in C_{c}^{\infty}({\mathbb{R}})$ satisfying  $(1- \varphi_{ij}^{\#}(x)) \varphi_{ij}(x)=0$ with $\varphi_{ij}(x)\equiv 1$ on a neighborhood of $x_0$. 

Here we recall $S^{m}$ the space of classical symbols of order $m$. Assuming $a\in S^{m}$, we denote $\op_{\e}(a)$ the pseudo-differential operators in semi-classical quantization who act on functions or distributions $u(x)$ with
\be
\op_{\e}(a)u=\int e^{ix\cdot \xi}a(x,\e\xi)\hat{u}(\xi)d\xi, \quad\e>0,\nn
\ee
and semi-classical norm $\|\cdot\|_{\e,s}$ with 
\be
\|u\|_{\e,s}^{2}=\int (1+|\e\xi|^{2})^{s}|\hat u(\xi)|^{2}d\xi.\nn
\ee
 It is rather clear that  there holds the identity $\op_{\e}(\sigma)(e^{ip\th}v)=e^{ip\th}\op_{\e}(\sigma_{+p})(v)$ with $\sigma_{+p}(x,\xi):=\sigma(x,\xi+pk)$.

For the case of coupled non-transparent resonances, we recall the following separation conditions.
\begin{assumption}\label{Separation conditions} Given $(i,j)\in \Re\backslash \Re_{0}$, the $(i,j)$-resonance is transparent. Given $(i,j)\in \Re_{0}$, the $(i,j)$-interaction coefficients are transparent on a neighborhood of
\be
{\mathcal{R}}_{ij}\bigcap(({\mathcal{R}}_{i'i}-k)\bigcup({\mathcal{R}}_{jj'}+k))\nn
\ee
for all $i',j'$ with $(i',i)\in \Re_{0}$, $(j,j')\in \Re_{0}$, and  on a neighborhood of
\be
{\mathcal{R}}_{ij}\bigcap({\mathcal{R}}_{ii'}\bigcup{\mathcal{R}}_{j'j})\nn
\ee
for all $i'\neq j, j'\neq i$, with $(i,i')\in \Re_{0}$, $(j',j)\in \Re_{0}$.
\end{assumption}
\begin{assumption}\label{Boundedness and Rank-one coefficients} Suppose 
\begin{itemize}
\item(Boundedness) The resonant set in ${\Re}$ is bounded.
\item(Rank-one coefficients) For all $(i,j)\in \Re_{0}$, for all $\xi$ in an open set containing ${\mathcal{R}}_{ij}$, the ranks of the $(i,j)$-interaction coefficients are at most 1; except for  $(i,j)\in \Re_{0}$ such that one interaction coefficient is identically equal to zero, in which case we make no assumptions on the rank of the other coefficient.
\end{itemize}
\end{assumption}
From the form of the stability index $\Gamma $ in \eqref{(1.11)}, it is much more likely that $\Gamma =1$ which corresponds to the instability result. Here we will focus on the instability case and  we recall Theorem 2.13 in \cite{Lu} concerning the instability of the WKB solutions.

In the context of Assumption \ref{Separation conditions}, for $(i,j)$ ranging over the set $\Re_{0}$ of non-transparent resonances, the maximum $\gamma$ of coefficients $\gamma_{ij}$ is attained at $(i_0,j_0)$. Consider the datum as 
\be\label{(1.15)}
u(0,x):=u_a(0,x)+\e^{k}e^{ix\cdot(\xi_0+k)/\e}\varphi_{i_0j_0}(x)\vec{e}_{i_0j_0},
\ee
where \par
--- $\xi_{0}$ is such that $\gamma:=|\max_{\xi\in{\mathcal{R}}_{i_0j_0}}\mbox{Re}(\Gamma_{i_0j_0}(\xi)^{\frac{1}{2}})|$ is attained at $\xi_0$;\par
--- $x_{0}$ is such that $|a|_{L^{\infty}}$ is attained at $x_{0}$;\par
--- $\varphi_{i_0j_0}\in C_{c}^{\infty}(\mathbb R^{d})$ is a scalar spatial truncation on a neighborhood of $x_0$;\par
--- $\vec{e}_{i_0j_0}$  is some fixed constant eigenvector corresponding to the eigenvalue of $M$ which generates exponential grow $O(e^{{t\gamma}^{+}})$.
\begin{theorem}\label{A instability criterion}Let $T_{\infty}:=\frac{K}{\gamma|a|_{\infty}}$, under Assumptions 1.1, 1.2, 1.4 and 1.5, in the unstable case ${\Gamma}>0$, for any $K>0$, if $K_{a}+\frac{1}{2}\geq K$, for $d/2<s\leq s_a$ :
\begin{itemize}
\item either for some $T < T_{\infty}$, for any $\e$ small enough, the initial-value problem \eqref{(1.1)}-\eqref{(1.15)} does not have a solution $u\in C^{0}([0,T\sqrt{\e}|\ln\e|], H^{s}({\mathbb{R}}^{d}))$;
\item or for some $T < T_{\infty}$, for any $\e_{0}>0$, the solution u to \eqref{(1.1)}-\eqref{(1.15)} satisfies
\be
\sup_{0<\e<\e_{0}}\sup_{0\leq t\leq T\sqrt{\e}|\ln\e|}|u(t,\cdot)|_{L^{\infty}}=\infty;\nn
\ee
\item or for any $K^{'}>0$, for some $T < T_{\infty}$, there holds the deviation estimate
\be
\sup_{0<\e<\e_{0}}\sup_{0\leq t\leq T\sqrt{\e}|\ln\e|}\e^{-K^{'}}|(u-u_{a})(t)|_{L^{2}}=\infty, \nn
\ee
\end{itemize}
for some $x_0\in \mathbb R^{d}$, some $\beta>0$, some $\e_{0}>0$.
\end{theorem}

\subsection{Estimates of the symbolic flow}
\label{Existence of the symbolic flow $2*2$}
The proof of the instability in \cite{Lu} of the WKB solutions when $\Gamma>0$ relies on a short-time Duhamel representation formula for solutions of zeroth-order pseudo-differential equations. To achieve such Duhamel representation formula, it is crucial to obtain the optimal upper bounds for the flows of the symbols of the pseudo-differential operators.  After frequency localizations  near $(i,j)$-resonance set and normal form reductions,  the separation conditions in Assumption \ref{Separation conditions} ensure that except the  $(i,j)$-interaction coefficients, all others will be eliminated. As a result,  the study of the symbolic flow reduces to the  local equations in $2 \times 2$ block matrices such as the following $M_{ij}$ in \eqref{(1.18)}.  We recall  the following lemma in \cite{Lu} which concludes the main results about the upper bounds for the symbolic flows of $2\times 2$ block matrices.
\begin{lemma}\label{Upper bounds for the symbolic flow $2*2$} Under Assumptions \ref{Smooth spectral decomposition}, \ref{WKB solution}, \ref{Separation conditions} and  \ref{Boundedness and Rank-one coefficients}, for all $T>0$ and all $ 0\leq\tau\leq t\leq T|\ln\e|$, $\alpha\in {\mathbb{N}}^{d}$, the solution $S_{0}$ to
 \be\label{(1.16)}
\partial_{t}S_{0}+\frac{1}{\sqrt{\e}}M_{ij}S_{0}=0, \quad S_{0}(\tau,\tau)=\Id
\ee
satisfies the bound
\be\label{(1.17)}
|\partial_{x}^{\alpha}S_{0}(\tau,t)|\lesssim|\ln\e|^{*}\exp(t\gamma_{ij}^{+})
\ee
with the block matrices of non-transparent coefficients $\{M_{ij}\}_{(i,j)\in \Re_{0}}$ defined as 
\ba\label{(1.18)}
M_{ij}(\e,t,x,\xi)={\chi^{\#}_{ij}}\bp  i(\l_{i,+1}-\omega) &  -\sqrt{\e}\varphi_{ij}^{\#}\chi_{ij}b_{ij}^{+}
    \\ 
    -\sqrt{\e}{\varphi^{\#}_{ij}}\chi_{ij}b_{ji}^{-} &  i\l_{j}
     \ep,
\ea
where 
\be\label{(1.19)}
b_{ij}^{+}=g(\sqrt{\e}t,x)\Pi_{i,+1}B_{1}\Pi_{j}, \quad b_{ji}^{-}=g^{*}(\sqrt{\e}t,x)\Pi_{j}B_{-1}\Pi_{i,+1}
\ee
and $1\leq i, j\leq J$, $i\neq j$. If the $M_{ij}$ in \eqref{(1.18)} is replaced by 
 \ba
M_{ij}(\e,t,x,\xi)={\chi^{\#}_{ij}}\bp  i\l_{i,+1} &  -\sqrt{\e}\varphi_{ij}^{\#}\chi_{ij}b_{ij}^{+}
    \\ 
    -\sqrt{\e}{\varphi^{\#}_{ij}}\chi_{ij}b_{ji}^{-} &  i(\l_{j}+\omega)\nn
     \ep,
\ea
the proof in \cite{Lu} shows the same bound for symbolic flows of  $2\times 2$ block matrices in \eqref{(1.17)}.  Furthermore, $|\ln\e|^{*}$ denotes $|\ln\e|^{N^*}$ for some $N^{*}>0$ independent of $(\e,\tau,t,x,\xi)$ and the support of $\chi_{ij}^{\#}$ is a small neighborhood of $\mathcal R_{ij}$. We specifically suppose $\chi_{ij}^{\#}$ be a smooth cut-off function in frequency space with $0\leq \chi_{ij}^{\#}\leq 1$, $\chi_{ij}^{\#}=1$ on a neighborhood of $\mathcal R_{ij}^{h}$ and $\chi_{ij}^{\#}=0$ away from $\mathcal R_{ij}^{h}$, such as $\mathbb R^{d}\backslash\mathcal R_{ij}^{2h}$.
 \end{lemma}
Assumption \ref{Separation conditions} concerning the separation conditions of the resonance sets is the key to reduce the massive matrix of interaction coefficients into some $2\times 2$ block matrices. After different normal form reduction and space-frequency localization, we now consider the following $3\times 3$ matrices cases   
 \ba\label{(1.20)}
M&=\bp  i\mu_{1} &  -\sqrt{\e}b_{12}  & -\sqrt{\e}b_{13}
    \\ 
     -\sqrt{\e}b_{21} & i\mu_{2} & -\sqrt{\e}b_{23}
    \\
    -\sqrt{\e}b_{31}  & -\sqrt{\e}b_{32} &  i\mu_{3} \ep.
\ea
Here we can not deal with the cases where all the interaction coefficients in $3\times 3$ block matrices  are non-transparent, while if some of the interaction coefficients in $3\times 3$ block matrices are transparent, the estimate of the symbolic flow can be obtained. Furthermore, we give relaxed  separation conditions and reasonable conditions in Assumption \ref{Reasonable conditions} to show those specific $3\times 3$ block matrices for which we can also obtain the estimate of symbolic flow finally.

\subsection{Main results}
\label{Main results}
For the sake of simplicity, we give the following reasonable conditions.
\begin{assumption}\label{Reasonable conditions} (Reasonable conditions) Assume $\l_{1}\geq\l_{2}\cdots\geq\l_{J}$, $\mathcal {R}_{ij}\bigcap(\mathcal {R}_{ij}+k)= \varnothing$ and $\forall \ (i,j)\in\Re_{0}\Rightarrow i<j$.
\end{assumption}

 Combining the property that the nonzero eigenvalue $\l_{i}(\xi)\rightarrow\infty$ as $|\xi|\rightarrow\infty$, we acquiescently choose $\l_{1}\geq\l_{2}\cdots\geq\l_{J}$ in this paper. We select a characteristic temporal frequency $\o\in \mathbb R$ associated with initial wavenumber $k\in \mathbb R^{d}$ satisfying that the phase $\beta=(\o,k)=(\l_{1}(k),k)$ belongs to the fastest positive branch on the variety. Then it is reasonable to assume the resonance occurs only when $i<j$ for $\mathcal R_{ij}$ (i.e. $\l_{i}\geq\l_{j}$ for $\mathcal R_{ij}$).  Thus under reasonable conditions in Assumption \ref{Reasonable conditions}, we can conclude following three cases for 2-coupled non-transparent resonances: $(i,j), (j,j^{'})$-resonances,  $(i,j), (i,j^{'})$-resonances,  and $(i,j^{'}), (j,j^{'})$-resonances, where $i<j<j^{'}$ and $\{(i,j), (i,j^{'}), (j,j^{'})\}\subset\Re_{0}$. Combined with these considerations,  we can reduce Assumption \ref{Separation conditions} to the following relaxed  separation conditions Assumption \ref{Relaxed  separation conditions}.
\begin{assumption}\label{Relaxed  separation conditions} (Relaxed  separation conditions)  Let $\Re_{0}\subset \Re$ be the collection of non-transparent pairs, i.e.  if $ (i,j)\in \Re\backslash \Re_{0}$, $(i,j)$-resonance is transparent. Assume that at least one of the following statements holds:
\begin{itemize}
\item[(1).] For any $i< j< j^{'}$ such that $(i,j)$, $(j,j^{'})\in \Re_{0}$, at least one of $b_{jj^{'}}^{+}$, $b_{j^{'}j}^{-}$ is transparent on $\mathcal {R}_{jj^{'}}\bigcap(\mathcal {R}_{ij}-k)$, or at least one of $b_{ij}^{+}$, $b_{ji}^{-}$ is transparent on $\mathcal {R}_{ij}\bigcap(\mathcal {R}_{jj^{'}}+k)$;
\item[(2).] For any $i< j< j^{'}$ such that  $(i,j)$, $(i,j^{'})\in \Re_{0}$, at least one of  $b_{ij^{'}}^{+}$, $b_{j^{'}i}^{-}$ is transparent on $\mathcal {R}_{ij}\bigcap\mathcal {R}_{ij^{'}}$, or at least one of  $b_{ij}^{+}$, $b_{ji}^{-}$ is transparent on $\mathcal {R}_{ij}\bigcap\mathcal {R}_{ij^{'}}$;
\item[(3).] For any $i< j< j^{'}$ such that  $(i,j^{'})$, $(j,j^{'})\in \Re_{0}$, at least one of $b_{jj^{'}}^{+}$, $b_{j^{'}j}^{-}$ is transparent on $\mathcal {R}_{ij^{'}}\bigcap\mathcal {R}_{jj^{'}}$, or at least one of $b_{ij^{'}}^{+}$, $b_{j^{'}i}^{-}$ is transparent on $\mathcal {R}_{ij^{'}}\bigcap\mathcal {R}_{jj^{'}}$.
\end{itemize}
\end{assumption}

Furthermore, we considered the 3-non-transparent resonances coupled in pairs which can be concluded as $\Re_{1}=\{(i,j),(i,j^{'}),(j,j^{'})\}\subset\Re_{0}$ and give the following assumption.
\begin{assumption}\label{3 coupled non-transparent resonances assumption}For non-transparent resonances in $\Re_{1}$, anyone of the following statements (separation conditions) holds: 
\begin{itemize}
\item[(1).] If $\mathcal {R}_{ij}\bigcap\mathcal {R}_{ij^{'}}\bigcap\mathcal {R}_{jj^{'}}\neq \varnothing$, $b_{jj^{'}}^{+}$ (or $b_{ji}^{-}$) is transparent in $\mathcal {R}_{jj^{'}}$ (or $\mathcal {R}_{ij}$) and anyone of  interaction coefficients in $\Re_{1}\backslash (j,j^{'})$ (or $\Re_{1}\backslash(i,j)$) is transparent on corresponding resonance set; 
\item[(2).] If $\mathcal {R}_{ij}\bigcap\mathcal {R}_{ij^{'}}\bigcap(\mathcal {R}_{jj^{'}}+k)\neq \varnothing$, $b_{j^{'}j}^{-}$ (or $b_{j^{'}i}^{-}$) is transparent in $\mathcal {R}_{jj^{'}}$ (or $\mathcal {R}_{ij^{'}}$) and anyone of  interaction coefficients in $\Re_{1}\backslash  (j,j^{'})$ (or $\Re_{1}\backslash (i,j^{'})$) is transparent on corresponding resonance set;
\item[(3).] If $\mathcal {R}_{ij}\bigcap(\mathcal {R}_{ij^{'}}+k)\bigcap(\mathcal {R}_{jj^{'}}+k)\neq \varnothing$, $b_{ij^{'}}^{+}$ (or $b_{ij}^{+}$) is transparent in $\mathcal {R}_{ij^{'}}$ (or $\mathcal {R}_{ij}$) and anyone of  interaction coefficients in $\Re_{1}\backslash  (i,j^{'})$ (or $\Re_{1}\backslash(i,j)$) is transparent on corresponding resonance set.
  \end{itemize}
\end{assumption}
Usually, we call $\mathcal R_{ij}$ the corresponding resonance set of interaction coefficients $b_{ij}^{+}$ and $b_{ji}^{-}$. For the cases of 2-coupled non-transparent resonances and 3-non-transparent resonances coupled in pairs, we have following new instability criterion for WKB solution $u_{a}$.
\begin{theorem}\label{New instability criterion}
Under reasonable conditions in Assumption \ref{Reasonable conditions}, we have the following results:
\begin{itemize}
\item[(i).] Theorem \ref{A instability criterion} holds if Assumption \ref{Separation conditions} in Theorem \ref{A instability criterion} is substituted by Assumption \ref{Relaxed  separation conditions};
\item[(ii).] Theorem \ref{A instability criterion} holds if Assumption \ref{Separation conditions} in Theorem \ref{A instability criterion} is substituted by Assumption \ref{3 coupled non-transparent resonances assumption}.
\end{itemize}
\end{theorem}
  
\subsection{Ideas of the proof}\label{Ideas of the proof}
This paper generalizes the main results in \cite{Lu} and allows non-separation of resonance sets. To overcome extra difficulties caused by the intersection of resonance sets,  there are two main new ideas compared to \cite{Lu}:
  
  \medskip
  
  1. Careful choices of frequency shit are given to delete oscillations $e^{ip\th}$ in $M_{ijj'}$. Since the oscillations $e^{ip\th}$ with interaction coefficients appear in the perturbed system, appropriate frequency shift is needed here to eliminate the oscillations. For example,  for any $i< j< j^{'}$ such that $(i,j)$, $(j,j^{'})\in \Re_{0}$, we give the frequency shift on $U_{i}$, $U_{j}$ and $U_{j^{'}}$ as $U_{i}=e^{-i\th}\op_{\e}{(\Pi_i)}\dot u$, $U_{j}=\op_{\e}{(\Pi_j)}\dot u$ and $ U_{j^{'}}=e^{i\th}\op_{\e}{(\Pi_{j^{'}})}\dot u$. Then the  perturbed system \eqref{(1.8)} can be diagonalized and furthermore, using normal form reduction and under  a $O(\sqrt{\e})$ remainder, the oscillations $e^{ip\th}$ are eliminated with interaction coefficients left supported in corresponding resonance sets. For this case, any other frequency shits are not ready to delete oscillations $e^{ip\th}$, that is the reason why we need to choose suitable frequency shit carefully. While for 3-non-transparent resonances coupled in pairs, frequency shift and normal form reduction can not delete all the oscillations. Under such conditions, we give several reasonable conditions to help remove the oscillations and finally it is reduced to one of 2-coupled non-transparent resonances. What's more, we find that different frequency shifts used in 2-coupled non-transparent resonances are also applicable for 3-non-transparent resonances coupled in pairs if we similarly change separation conditions, which offer us two more possibilities to give the final instability of WKB approximate solution $u_{a}$. 
  
  \medskip

  2. Estimates for the symbolic flow related to $3\times 3$ block matrices $M_{ijj'}$ are established. After normal form reduction, non-transparent interaction coefficients supported in corresponding resonance sets are left without oscillations, which provides us the convenience for the following space-frequency localization. Utilizing the estimates for symbolic flow related to $2\times 2$ block matrices $\{M_{ij}\}_{(i,j)\in \Re_{0}}$, which is given in Lemma \ref{Upper bounds for the symbolic flow $2*2$}, the key point is the classified discussion in different non-transparent sets and their intersections. Especially in their intersections, separation conditions will help reduce one of the non-transparent interaction coefficients and  successfully give the estimates for the symbolic flow. Whether in cases of 2-coupled non-transparent resonances or 3-non-transparent resonances coupled in pairs, classified discussion combined with separation conditions and reasonable conditions show the estimates for the symbolic flow and finally give the instability by the same method of Duhamel representation in \cite{Lu}.
  
  \begin{remark}\label{Remark on 3-non-transparent resonances}  
For 3-non-transparent resonances coupled in pairs in $\Re_{1}$ in Assumption \ref{3 coupled non-transparent resonances assumption}, we choose suitable reasonable condition $\mathcal {R}_{ij}\bigcap(\mathcal {R}_{ij}+k)=\varnothing$ and separation conditions to wipe out the oscillation items and reduce the case to one of the cases talked above in Assumption \ref{Relaxed  separation conditions} describing 2-coupled non-transparent resonances. From Theorem \ref{A instability criterion} and Theorem \ref{New instability criterion}, it can be observed that if non-transparent resonance set $\Re_{0}$ for the highly-oscillating solutions to semi-linear systems is composed of  $\Re_{1}$ satisfying conditions in Assumption \ref{3 coupled non-transparent resonances assumption}, 2-coupled non-transparent resonances satisfying relaxed  separation conditions in Assumption \ref{Relaxed  separation conditions} and 1-single non-transparent resonance, we can also verify the WKB approximate solution $u_{a}$ satisfies the instability result in Theorem \ref{A instability criterion}.
\end{remark}

\section{Cases of two coupled resonances}
\label{Cases of two coupled resonances}
  To prove the first result in Theorem \ref{New instability criterion}, there is no influence if we consider the following specific assumptions of separation conditions. Given $(i,j,j^{'})=(1,2,3)$,  Assumption \ref{Relaxed  separation conditions} becomes:
\begin{itemize}
\item[(1).] Assume $(i,j)=(1,2)$, $(j,j^{'})=(2,3)\in \Re_{0}$, at least one of $b_{23}^{+}$, $b_{32}^{-}$ is transparent on $\mathcal {R}_{23}\bigcap(\mathcal {R}_{12}-k)$, or at least one of $b_{12}^{+}$, $b_{21}^{-}$ is transparent on $\mathcal {R}_{12}\bigcap(\mathcal {R}_{23}+k)$;
\item[(2).] Assume $(i,j)=(1,2)$, $(i,j^{'})= (1,3)\in \Re_{0}$, at least one of  $b_{13}^{+}$, $b_{31}^{-}$ is transparent on $\mathcal {R}_{12}\bigcap\mathcal {R}_{13}$, or at least one of  $b_{12}^{+}$, $b_{21}^{-}$ is transparent on $\mathcal {R}_{12}\bigcap\mathcal {R}_{13}$;
\item[(3).] Assume $(i,j^{'})=(1,3)$, $(j,j^{'})=(2,3)\in \Re_{0}$, at least one of $b_{23}^{+}$, $b_{32}^{-}$ is transparent on $\mathcal {R}_{13}\bigcap\mathcal {R}_{23}$, or at least one of $b_{13}^{+}$, $b_{31}^{-}$ is transparent on $\mathcal {R}_{13}\bigcap\mathcal {R}_{23}$.
\end{itemize}

Now we need to show the instability in Theorem \ref{A instability criterion} holds under Assumption \ref{Smooth spectral decomposition}, Assumption \ref{WKB solution}, Assumption \ref{Boundedness and Rank-one coefficients} and furthermore, any one of above three cases of assumptions.
\subsection{Case 1}\label{Case 1}
  We start by considering the first case  where $(1,2)$, $(2,3)\in \Re_{0}$, at most one of $b_{23}^{+}$, $b_{32}^{-}$ is non-transparent on $\mathcal {R}_{23}\bigcap(\mathcal {R}_{12}-k)$, or at most one of $b_{12}^{+}$, $b_{21}^{-}$ is non-transparent on $\mathcal {R}_{12}\bigcap(\mathcal {R}_{23}+k)$.
  
\subsubsection{Diagonalizaiton and shift of frequency}

By the eigenmodes of the hyperbolic operator, we decompose $\dot u$ and shift the component related to $\Pi_{1}$ and  $\Pi_{3}$. Then we make the following frequency shift
\be
U_{1}=e^{-i\th}\op_{\e}{(\Pi_1)}\dot u,  U_{2}=\op_{\e}{(\Pi_2)}\dot u,  U_{3}=e^{i\th}\op_{\e}{(\Pi_3)}\dot u,U_{4}=\op_{\e}{(\Pi_4)}\dot u,\cdots, U_{J}=\op_{\e}{(\Pi_J)}\dot u,\nn
\ee
where the projectors $\Pi_{j} (j=1,2,3,\cdots,J)$ are eigenprojectors of $A(i\xi)+A_0$ and $op_\e(\Pi_j) (j=1,2,3,\cdots,J)$ are defined as above. The perturbation unknown $\dot u$ is decomposed by
\be
\dot u=e^{i\th}U_{1}+U_{2}+e^{-i\th}U_{3}+U_{4}+\cdots+U_{J}.\nn
\ee
Then we find that $U=(U_1,U_2,U_3,\cdots,U_J)$ solves
\be\label{(2.1)}
\partial_{t}U+\frac{1}{\e}\op_{\e}(i\mathcal A)U=\frac{1}{\sqrt{\e}}\op_{\e}(i\mathcal B)U+F.
\ee
The symbol of the propagator is 
\be\label{(2.2)}
\mathcal A=\diag(\l_{1,+1}-w,\l_{2},\l_{3,-1}+w,\l_{4},\cdots,\l_{J}),
\ee
and the symbol of singular source term is

     \ba\label{(2.3)}
     \mathcal B=\bp  \mathcal B_{[1,3]} &\mathcal B_{[1,3,J]}\\ 
     \mathcal B_{[J,1,3]} & \mathcal B_{[J,J]}\ep
     \ea
     where the top left block is
     \ba
  \mathcal B_{[1,3]}&=\sum_{p=\pm 1}
  \bp  e^{ip\th}\Pi_{1,+(p+1)}B_{p}\Pi_{1,+1} &  e^{i(p-1)\th}\Pi_{1,+p}B_{p}\Pi_{2}  &  e^{i(p-2)\th}\Pi_{1,+(p-1)}B_{p}\Pi_{3,-1}
    \\ 
     e^{i(p+1)\th}\Pi_{2,+(p+1)}B_{p}\Pi_{1,+1} & e^{ip\th}\Pi_{2,+p}B_{p}\Pi_{2} & e^{i(p-1)\th}\Pi_{2,+(p-1)}B_{p}\Pi_{3,-1}
    \\
     e^{i(p+2)\th}\Pi_{3,+(p+1)}B_{p}\Pi_{1,+1} & e^{i(p+1)\th}\Pi_{3,+p}B_{p}\Pi_{2} &  e^{ip\th}\Pi_{3,+(p-1)}B_{p}\Pi_{3,-1}\ep\in \mathbb R^{3N\times 3N},\nn
     \ea
       with $\Pi_{j,+q}{(\xi)}=\Pi(\xi+qk)$ for $q\in\ZZ$ and $B_{p}=B(u_{0,p}), ~p\in{\{-1,1\}}$, where $u_{0,\pm 1}$ are the leading amplitudes in the WKB approximate solution. The other blocks are     
           \ba
  \mathcal B_{[1,3,J]}&=\sum_{p=\pm 1}
  \bp   e^{i(p-1)\th}\Pi_{1,+p}B_{p}\Pi_{4}&\cdots &  e^{i(p-1)\th}\Pi_{1,+p}B_{p}\Pi_{J}
    \\ 
      e^{ip\th}\Pi_{2,+p}B_{p}\Pi_{4}&\cdots &  e^{ip\th}\Pi_{2,+p}B_{p}\Pi_{J}
    \\
       e^{i(p+1)\th}\Pi_{3,+p}B_{p}\Pi_{4}&\cdots &  e^{i(p+1)\th}\Pi_{3,+p}B_{p}\Pi_{J}\ep\in \mathbb R^{3N\times (J-3)N},\nn
     \ea
\ba
  \mathcal B_{[J,1,3]}&=\sum_{p=\pm 1}
  \bp  e^{i(p+1)\th}\Pi_{4,+(p+1)}B_{p}\Pi_{1,+1} &  e^{ip\th}\Pi_{4,+(p-1)}B_{p}\Pi_{2,-1}  &  e^{i(p-1)\th}\Pi_{4,+(p-1)}B_{p}\Pi_{3,-1}
     \\ 
     \vdots & \vdots & \vdots
     \\
     e^{i(p+1)\th}\Pi_{J,+(p+1)}B_{p}\Pi_{1,+1} &  e^{ip\th}\Pi_{J,+(p-1)}B_{p}\Pi_{2,-1}  &  e^{i(p-1)\th}\Pi_{J,+(p-1)}B_{p}\Pi_{3,-1}\ep\in \mathbb R^{(J-3)N\times 3N},\nn
     \ea
     and   
  \ba
  \mathcal B_{[J,J]}&=\sum_{p=\pm 1}
  \bp   e^{ip\th}\Pi_{4,+p}B_{p}\Pi_{4}&\cdots & e^{ip\th}\Pi_{4,+p}B_{p}\Pi_{J}
    \\ 
      \vdots& &  \vdots
    \\
       e^{ip\th}\Pi_{J,+p}B_{p}\Pi_{4}&\cdots &  e^{ip\th}\Pi_{J,+p}B_{p}\Pi_{J}\ep\in \mathbb R^{(J-3)N\times (J-3)N}.\nn
     \ea  
 We decompose the top left block of $\mathcal B$ as
 \be
 \mathcal B_{[1,3]}=\mathcal B^{r}+\mathcal B^{nr},\nn
 \ee
with the notation
   \ba\label{(2.4)}
  \mathcal B^{r}=\bp  0 &  \Pi_{1,+1}B_{1}\Pi_{2}  &  0
    \\ 
     \Pi_{2}B_{-1}\Pi_{1,+1} &0 & \Pi_{2}B_{1}\Pi_{3,-1}
    \\
    0 & \Pi_{3,-1}B_{-1}\Pi_{2} & 0\ep.
     \ea  
Furthermore, we denote  
     \ba\label{(2.5)}
      \mathcal D=\bp  \tilde{\mathcal B}+\mathcal B^{nr} &\mathcal B_{[1,3,J]}\\ 
     \mathcal B_{[J,1,3]} & \mathcal B_{[J,J]}\ep,
     \ea
with 
  \ba
  \tilde{\mathcal B}&=\bp  0 &  (1-\chi_{12})\Pi_{1,+1}B_{1}\Pi_{2}  &  0
    \\ 
     (1-\chi_{12})\Pi_{2}B_{-1}\Pi_{1,+1} &0 & (1-\chi_{23,-1})\Pi_{2}B_{1}\Pi_{3,-1}
    \\
    0 & (1-\chi_{23,-1})\Pi_{3,-1}B_{-1}\Pi_{2} & 0\ep.\nn
     \ea  
 Assumptions in Theorem \ref{A instability criterion} ensure the interaction coefficients in $\mathcal B^{nr}$, $\mathcal B_{[1,3,J]}$, $\mathcal B_{[J,1,3]}$ and $\mathcal B_{[J,J]}$ are transparent with oscillations and the interaction coefficients in $\tilde{\mathcal B}$ are supported away from corresponding resonance sets. In a skillful way, we will delete them using the following normal form reduction, which means we can delete $\mathcal D$ in the system with a $O(\sqrt{\e})$ remainder.
  \subsubsection{Normal form reduction} 
Using normal form reduction with $\mu_{1}=\l_{1,+1}-\omega, \ \mu_2=\l_2, \ \mu_{3}=\l_{3,-1}+\omega,\ \mu_j=\l_j$ for $4\leq j \leq J$, it is sufficient to solve 
\be\label{(2.6)}
i(-p\omega+\mu_{i,+p}-\mu_{j})({Q_{p}})_{(i,j)}=({\mathcal D_{p}})_{(i,j)}, \quad 1\leq i,j \leq J,
\ee
 where we rewrite the above $\mathcal D$ as $\mathcal D=\sum_{|p|\leq 3}e^{ip\th}{\mathcal D_{p}}$.
 
 When $p=0$, we have
 \be
i(\mu_1-\mu_2)(Q_0)_{(1, 2)}=(1-\chi_{12})\Pi_{1,+1}B_{1}\Pi_{2}\Rightarrow i(\l_{1,+1}-\l_2-\omega)(Q_0)_{(1, 2)}=(1-\chi_{12})b_{12}^{+};\nn
\ee
\be
i(\mu_2-\mu_1)(Q_0)_{(2, 1)}=(1-\chi_{12})\Pi_{2}B_{-1}\Pi_{1,+1}\Rightarrow -i(\l_{1,+1}-\l_2-\omega)(Q_0)_{(2, 1)}=(1-\chi_{12})b_{21}^{-};\nn
\ee
\be
i(\mu_2-\mu_3)(Q_0)_{(1, 2)}=(1-\chi_{23,-1})\Pi_{2}B_{1}\Pi_{3,-1}\Rightarrow i(\l_{2,+1}-\l_{3}-\omega)_{-1}(Q_0)_{(2, 3)}=(1-\chi_{23,-1})(b_{23}^{+})_{-1};\nn
\ee
\be
i(\mu_3-\mu_2)(Q_0)_{(1, 2)}=(1-\chi_{23,-1})\Pi_{3,-1}B_{-1}\Pi_{2}\Rightarrow -i(\l_{2,+1}-\l_{3}-\omega)_{-1}(Q_0)_{(3, 2)}=(1-\chi_{23,-1})(b_{32}^{-})_{-1};\nn
\ee
\be
i(\mu_1-\mu_j)(Q_0)_{(1, j)}=b_{1j}^{+}\Rightarrow i(\l_{1,+1}-\l_j-\omega)(Q_0)_{(1, j)}=b_{1j}^{+},\quad ( 4\leq j\leq J);\nn
\ee
\be
i(\mu_3-\mu_j)(Q_0)_{(3, j)}=(b_{3j}^{-})_{-1}\Rightarrow -i(\l_{j,+1}-\l_3-\omega)_{-1}(Q_0)_{(3, j)}=(b_{3j}^{-})_{-1},  \quad( 4\leq j\leq J);\nn
\ee
\be
i(\mu_j-\mu_1)(Q_0)_{(j, 1)}=b_{j1}^{-}\Rightarrow -i(\l_{1,+1}-\l_j-\omega)(Q_0)_{(j, 1)}=b_{j1}^{-},\quad  ( 4\leq j\leq J);\nn
\ee
\be
i(\mu_j-\mu_3)(Q_0)_{(j,3)}=(b_{j3}^{+})_{-1}\Rightarrow i(\l_{j,+1}-\l_3-\omega)_{-1}(Q_0)_{(j, 3)}=(b_{j3}^{+})_{-1}, \quad (4\leq j\leq J).\nn
\ee
In $\supp(1-\chi_{12})$,  $|(\l_{1,+1}-\l_2-\omega)|$ has a maximum positive lower bound and in $\supp(1-\chi_{23,-1})$, $|(\l_{2,+1}-\l_{3}-\omega)_{-1}|$ has a maximum positive lower bound too. Then we can divide the right-hand side by the phase to define $(Q_0)_{(1, 2)}$,  $(Q_0)_{(2, 1)}$, $(Q_0)_{(2, 3)}$ and  $(Q_0)_{(3, 2)}$ in the above cases. Relaxed  separation conditions in Assumption \ref{Relaxed  separation conditions} ensure the phases factor out in the right-hand sides, then we can solve $(Q_0)_{(1, j)}$, $(Q_0)_{(3, j)}$, $(Q_0)_{(j,1)}$ and $(Q_0)_{(j,3)}$ for $( 4\leq j\leq J)$ in $S^{0}$.

When $|p|=1$, we have
\be
i(\l_{1,+1}-\l_1-\omega)_{+1}(Q_1)_{(1, 1)}=\Pi_{1,+2}B_{1}\Pi_{1,+1}=(b_{11}^{+})_{+1},\nn
\ee
 \be
-i(\l_{1,+1}-\l_1-\omega)(Q_{-1})_{(1, 1)}=\Pi_{1}B_{-1}\Pi_{1,+1}=b_{11}^{-},\nn
\ee
 \be
i(\l_{1,+1}-\l_3-\omega)_{-1}(Q_{-1})_{(1, 3)}=\Pi_{1}B_{1}\Pi_{3,-1}=(b_{13}^{+})_{-1},\nn
\ee
 \be
i(\l_{2,+1}-\l_2-\omega)(Q_{1})_{(2, 2)}=\Pi_{2,+1}B_{1}\Pi_{2}=b_{22}^{+},\nn
\ee
 \be
-i(\l_{2,+1}-\l_2-\omega)_{-1}(Q_{-1})_{(2, 2)}=\Pi_{2,-1}B_{-1}\Pi_{2}=(b_{22}^{-})_{-1},\nn
\ee
 \be
-i(\l_{1,+1}-\l_3-\omega)(Q_{1})_{(3, 1)}=\Pi_{3}B_{-1}\Pi_{1,+1}=b_{31}^{-},\nn
\ee
 \be
i(\l_{3,+1}-\l_3-\omega)_{-1}(Q_{1})_{(3, 3)}=\Pi_{3}B_{1}\Pi_{3,-1}=(b_{33}^{+})_{-1},\nn
\ee
\be
-i(\l_{3,+1}-\l_3-\omega)_{-2}(Q_{-1})_{(3, 3)}=\Pi_{3,-2}B_{-1}\Pi_{3,-1}=(b_{33}^{-})_{-2},\nn
\ee
 \be
i(\l_{i,+1}-\l_j-\omega)(Q_{1})_{(i, j)}=\Pi_{i,+1}B_{1}\Pi_{j}=b_{ij}^{+}, \quad ( 4\leq i,j\leq J),\nn
\ee
\be
-i(\l_{j,+1}-\l_i-\omega)_{-1}(Q_{-1})_{(i, j)}=\Pi_{i,-1}B_{-1}\Pi_{j}=(b_{ij}^{-})_{-1}, \quad ( 4\leq i,j\leq J).\nn
\ee
When $|p|=2$, we have
\be
-i(\l_{2,+1}-\l_1-\omega)_{-1}(Q_{-2})_{(1, 2)}=\Pi_{1,-1}B_{-1}\Pi_{2}=(b_{12}^{-})_{-1},\nn
\ee
\be
i(\l_{2,+1}-\l_1-\omega)_{+1}(Q_{2})_{(2, 1)}=\Pi_{2,+2}B_{1}\Pi_{1+1}=(b_{21}^{+})_{+1},\nn
\ee
\be
-i(\l_{3,+1}-\l_2-\omega)_{-2}(Q_{-2})_{(2, 3)}=\Pi_{2,-2}B_{-1}\Pi_{3,-1}=(b_{23}^{-})_{-2},\nn
\ee
\be
i(\l_{3,+1}-\l_2-\omega)(Q_{2})_{(3, 2)}=\Pi_{3,-2}B_{1}\Pi_{2}=b_{32}^{+},\nn
\ee
\be
-i(\l_{j,+1}-\l_1-\omega)_{-1}(Q_{-2})_{(1, j)}=\Pi_{1,-1}B_{-1}\Pi_{j}=(b_{1j}^{-})_{-1}, \quad ( 4\leq j\leq J),\nn
\ee
\be
i(\l_{3,+1}-\l_j-\omega)_{-2}(Q_{-2})_{(j,3)}=\Pi_{j,-2}B_{-1}\Pi_{3,-1}=(b_{j3}^{-})_{-2}, \quad ( 4\leq j\leq J),\nn
\ee
 \be
i(\l_{3,+1}-\l_j-\omega)(Q_{2})_{(3, j)}=\Pi_{3,+1}B_{1}\Pi_{j}=b_{3j}^{+}, \quad ( 4\leq j\leq J),\nn
\ee
\be
i(\l_{j,+1}-\l_{1}-\omega)_{+1}(Q_{2})_{(j, 1)}=\Pi_{j,+2}B_{1}\Pi_{1,+1}=(b_{j1}^{+})_{+1}, \quad ( 4\leq j\leq J).\nn
\ee
When $|p|=3$, we have
\be
-i(\l_{3,+1}-\l_1-\omega)_{-2}(Q_{-3})_{(1, 3)}=\Pi_{1,-2}B_{-1}\Pi_{3,-1}=(b_{13}^{-})_{-2},\nn
\ee
\be
i(\l_{3,+1}-\l_{1}-\omega)_{+1}(Q_{3})_{(3, 1)}=\Pi_{3,+2}B_{1}\Pi_{1,+1}=(b_{31}^{+})_{+1}.\nn
\ee
Relaxed  separation conditions in Assumption \ref{Relaxed  separation conditions} are used again to solve for $Q$  in $S^{0}$ in above cases of $|p|=1,2,3$, and combining them. We solve the equations \eqref{(2.6)} and obtain $({Q_{p}})_{(i,j)}\in S^{0}$ for $1\leq i,j \leq J$ and $|p|=0,1,2,3$ in corresponding support sets. Given the change of variable $\check U(t)=(\Id+\sqrt{\e}\op_{\e}(Q(\sqrt{\e}t)))^{-1}U(\sqrt{\e}t)$ and combining the above normal form reduction similar to Corollary 3.5 of \cite{Lu}, we find equation in $\check U(t)$ :
 \be\label{(2.7)}
\partial_{t}\check U+\frac{1}{\sqrt{\e}}\op_{\e}(i\mathcal A)\check U=\op_{\e}(\check{\mathcal B})\check U+\sqrt{\e}\check F,
\ee
where the symbol $\check{B}$ is as 
\ba\label{(2.8)}
  \check{B}&=\bp  0 &  \chi_{12}\Pi_{1,+1}B_{1}\Pi_{2}  & 0
    \\ 
     \chi_{12}\Pi_{2}B_{-1}\Pi_{1,+1} & 0 & \chi_{23,-1}\Pi_{2}B_{1}\Pi_{3,-1}
    \\
    0   & \chi_{23,-1}\Pi_{3,-1}B_{-1}\Pi_{2} &  0 \ep\\
   & :=\bp  0 &  \chi_{12}b_{12}^{+}  & 0
    \\ 
     \chi_{12}b_{21}^{-} & 0 & \chi_{23,-1}(b_{23}^{+})_{-1}
    \\
    0   & \chi_{23,-1}(b_{32}^{-})_{-1} &  0 \ep.
  \ea  
  \subsubsection{Space-frequency localization}
Furthermore, we let
\be\label{(2.9)}
V=\op_{\e}(\chi)(\varphi\check{U}),\quad W=(W_{1},W_{2})(\op_{\e}(\chi)((1-\varphi) \check{U}), (1-\op_{\e}(\chi))\check{U})
\ee
satisfying
\be
\check{U}=V+W_{1}+W_{2}.\nn
\ee

Then we have
  \be\label{(2.10)}
\partial_{t}V+\frac{1}{\sqrt{\e}}\op_{\e}^{\psi}(M)V=\sqrt{\e}F_{V},
\ee
\be\label{(2.11)}
\partial_{t}W+\frac{1}{\sqrt{\e}}\op_{\e}(iA)W=\op_{\e}(D)W+\sqrt{\e}F_{W}
\ee
with symbol
\be\label{(2.12)}
M=i\chi^{\#} \mathcal A-\sqrt{\e}\varphi^{\#}  \check{B}, \quad
 A=\bp \mathcal{A} & 0
 \\
 0 &\mathcal{A}\ep, \quad
 D=(1-\varphi)\chi^{\#}\check{B}.
 \ee
$\check F$, $F_{W}$ and $F_{V}$ satisfy the same estimates in Lemma 3.1 of \cite{Lu}. While in different frequency sets $\supp\chi(\xi)$ and spatial truncations $\supp\varphi(x)$, small differences in the symbols $M, \check{B}$ and $(V,W)$ occur, which will be specified in following discussions. 
\subsubsection{Estimates of symbolic flow}
Now, we are trying to obtain the bounds for symbolic flow in the following equations
\be\label{(2.13)}
\partial_{t}S_{0}+\frac{1}{\sqrt{\e}}MS_{0}=0, \quad S_{0}(\tau,\tau)=\Id.
\ee
In other words, we need to give estimate for $|S_{0}(\tau,t)|=|\exp(-\frac{t}{\sqrt{\e}}M)|$ and most importantly, we want to obtain 
\be\label{(2.14)}
|S_{0}(\tau,t)|\lesssim|\ln\e|^{*}\exp(t\gamma^{+})
\ee
 for frequency $\xi$ in any non-transparent resonance sets or not. As discussed in \cite{Lu}, if $\xi$ is in transparent resonance sets or no resonance occurs, the normal form reduction will eliminate corresponding effects in \eqref{(2.10)} and \eqref{(2.11)} with two small remainders on the right side. Hence we need to discuss when  frequency $\xi$ in different non-transparent sets, under certain separation conditions (i.e. separation conditions in Assumption \ref{Relaxed  separation conditions}), whether the bound for symbolic flow in $\eqref{(2.14)}$ holds or not.

 If $\mathcal{R}_{12}\cap(\mathcal{R}_{23}+k)=\varnothing$ and $\xi\in\mathcal{R}_{12}$, let $\chi=\chi_{12}$ and $\varphi=\varphi_{12}$ in \eqref{(2.9)} and \eqref{(2.12)}, then we have
 \ba
  \check{B}=\chi_{12}\bp  0 & b_{12}^{+}  & 0
    \\ 
     b_{21}^{-} & 0 & 0
    \\
    0   & 0 &  0 \ep,\nn
  \ea
  and
  \ba\label{(2.15)}
M&=i\chi_{12}^{\#} \mathcal A-\sqrt{\e}\varphi_{12}^{\#}  \check{B}\\
 &=\bp  i\chi_{12}^{\#}(\l_{1,+1}-\omega) &  -\sqrt{\e}\varphi_{12}^{\#}\chi_{12}b_{12}^{+}  & 0
    \\ 
     -\sqrt{\e}\varphi_{12}^{\#}\chi_{12}b_{21}^{-} & i\chi_{12}^{\#}\l_{2} & 0
    \\
    0   & 0 &  i\chi_{12}^{\#}(\l_{3,-1}+\omega) \ep,
\ea
  then \eqref{(2.13)} can be solved directly using Lemma \ref{Upper bounds for the symbolic flow $2*2$} with $|S_{0}(\tau,t)|\lesssim|\ln\e|^{*}\exp(t\gamma_{12}^{+})\lesssim|\ln\e|^{*}\exp(t\gamma^{+})$ for $\xi\in \mathcal{R}_{12}$. 
  
  And for $\xi\in(\mathcal{R}_{23}+k)$, let $\chi=\chi_{23,-1}$, $\varphi=\varphi_{23}$ in \eqref{(2.9)} and \eqref{(2.12)}, then we have
  \ba\label{(2.16)}
  \check{B}=\chi_{23,-1}\bp  0 &  0  & 0
    \\ 
     0 & 0 & (b_{23}^{+})_{-1}
    \\
    0   & (b_{32}^{-})_{-1} &  0 \ep
   :=\chi_{23}{(\xi-k)}\bp  0 &  0  & 0
    \\ 
     0 & 0 & b_{23}^{+}{(\xi-k)}
         \\
    0   & b_{32}^{-}{(\xi-k)} &  0 \ep.
  \ea
 \\
  Let $\xi^{'}=\xi-k$, then $\xi^{'}\in\mathcal{R}_{23}$ and we have
  \ba
    \check{B}=\chi_{23}{(\xi^{'})}\bp  0 &  0  & 0
    \\ 
     0 & 0 & b_{23}^{+}{(\xi^{'})}
    \\
    0   & b_{32}^{-}{(\xi^{'})} &  0 \ep. \nn
  \ea
For the equation
\be
\partial_{t}S_{0}+\frac{1}{\sqrt{\e}}M(\xi^{'})S_{0}=0, \quad S_{0}(\tau,\tau)=\Id\nn
\ee
with
\ba
M(\xi^{'})&=i\chi^{\#}_{23,-1} A-\sqrt{\e}\varphi^{\#}_{23} \check{B}\\
&=i\chi^{\#}_{23}(\xi^{'})
\bp  \l_{1}(\xi^{'}+2k)-\omega &  0  & 0
    \\ 
     0 & \l_{2}(\xi^{'}+k) & 0
    \\
    0   & 0 &  \l_{3}(\xi^{'})+\omega \ep
    -\sqrt{\e}\varphi^{\#}_{23} \chi_{23}{(\xi^{'})}
    \bp  0 &  0  & 0
    \\ 
     0 & 0 & b_{23}^{+}{(\xi^{'})}
    \\
    0   & b_{32}^{-}{(\xi^{'})} &  0 \ep\\
    &=\bp  i\chi^{\#}_{23}(\xi^{'})(\l_{1}(\xi^{'}+2k)-\omega) &  0  & 0
    \\ 
     0 & i\chi^{\#}_{23}(\xi^{'})\l_{2}(\xi^{'}+k) & -\sqrt{\e}\varphi^{\#}_{23} \chi_{23}{(\xi^{'})}b_{23}^{+}{(\xi^{'})}
    \\
    0   & -\sqrt{\e}\varphi^{\#}_{23} \chi_{23}{(\xi^{'})}b_{23}^{+}{(\xi^{'})} &  i\chi^{\#}_{23}(\xi^{'})(\l_{3}(\xi^{'})+\omega) \ep,\nn
\ea
using the Lemma \ref{Upper bounds for the symbolic flow $2*2$}, we obtain
\be\label{(2.17)}
|S_{0}(\tau,t)|\lesssim|\ln\e|^{*}\exp(t\gamma_{23}^{+}(\xi^{'}))=|\ln\e|^{*}\exp(t\gamma_{23}^{+}).
\ee
where
\be
\gamma_{23}^{+}:=|a|_{L^{\infty}}|\max_{\xi\in{\mathcal{R}}_{23}^{h}}\mbox{Re}(\Gamma_{23}(\xi)^{\frac{1}{2}})|,\nn
\ee
\be
\gamma^{+}_{23}(\xi^{'})=|a|_{L^{\infty}}|\max_{\xi^{'}\in{\mathcal{R}}_{23}^{h}}\mbox{Re}(\Gamma_{23}(\xi^{'})^{\frac{1}{2}})|=|a|_{L^{\infty}}|\max_{\xi\in{\mathcal{R}}_{23}^{h}}\mbox{Re}(\Gamma_{23}(\xi)^{\frac{1}{2}})|=\gamma_{23}^{+}.\nn
\ee
Then by Lemma \ref{Upper bounds for the symbolic flow $2*2$}, we get the estimate $|S_{0}(\tau,t)|\lesssim|\ln\e|^{*}\exp(t\gamma_{23}^{+})\lesssim|\ln\e|^{*}\exp(t\gamma^{+})$ for $\xi\in (\mathcal{R}_{23}+k)$, where  $\mathcal{R}_{12}\cap(\mathcal{R}_{23}+k)=\varnothing$.

From \eqref{(2.16)} to \eqref{(2.17)}, the following Remark can be concluded.
\begin{remark}\label{Remark 2.1} Variable substitution in \eqref{(2.16)} to \eqref{(2.17)} shows the same upper bound for the symbolic flow. Similarly, we can conclude that if there is a uniform Variable substitution in frequency parameter i.e. for
\ba
M_{ij}(\e,t,x,\xi)={\chi^{\#}_{ij,-1}}\bp  i\l_{i} &  -\sqrt{\e}\varphi_{ij}^{\#}\chi_{ij}(b_{ij}^{+})_{-1}
    \\ 
    -\sqrt{\e}{\varphi^{\#}_{ij}}\chi_{ij}(b_{ji}^{-})_{-1} &  i(\l_{j,-1}+\omega)
     \ep,\nn
\ea
corresponding the phase to $\l_{i}-(\l_{j,-1}+\omega)=(\l_{i,+1}-\l_{j}+\omega)_{-1}$ and following the proof in \cite{Lu}, we can obtain \eqref{(2.14)} for \eqref{(2.13)} too.  
\end{remark}
  If $\mathcal{R}_{12}\cap(\mathcal{R}_{23}+k)\neq\varnothing$, for $\xi\in\mathcal{R}_{12}\backslash(\mathcal{R}_{23}+k)$, let $\chi=\chi_{{{\mathcal R}_{12}\backslash(\mathcal{R}_{23}+k)}}$ and $\varphi=\varphi_{12}$ in \eqref{(2.9)} and \eqref{(2.12)}, then we have
  \ba\label{(2.18)}
  \check{B}=\chi_{{\mathcal R}_{12}\backslash(\mathcal{R}_{23}+k)}\bp  0 & b_{12}^{+}  & 0
    \\ 
     b_{21}^{-} & 0 & 0
    \\
    0   & 0 &  0 \ep, 
  \ea
  and for $\xi\in(\mathcal{R}_{23}+k)\backslash\mathcal{R}_{12}$, let $\chi=\chi_{(\mathcal{R}_{23}+k)\backslash\mathcal{R}_{12}}$ and $\varphi=\varphi_{23}$ in \eqref{(2.9)} and \eqref{(2.12)}, then we have
   \ba
  \check{B}=\chi_{(\mathcal{R}_{23}+k)\backslash\mathcal{R}_{12}}
  \bp  0 &  0  & 0
    \\ 
     0 & 0 & (b_{23}^{+})_{-1}
    \\
    0   & (b_{32}^{-})_{-1} &  0 \ep
    =\chi_{(\mathcal{R}_{23+k})\backslash\mathcal{R}_{12}}(\xi)
  \bp  0 &  0  & 0
    \\ 
     0 & 0 & b_{23}^{+}(\xi-k)
    \\
    0   & b_{32}^{-}(\xi-k) &  0 \ep.\nn
  \ea
  Let $\xi^{'}=\xi-k$, then $\xi^{'}\in\mathcal{R}_{23}\backslash(\mathcal{R}_{12}-k)$ and we have
     \ba\label{(2.19)}
  \check{B}=\chi_{\mathcal{R}_{23}\backslash(\mathcal{R}_{12}-k)}(\xi^{'})
  \bp  0 &  0  & 0
    \\ 
     0 & 0 & b_{23}^{+}(\xi^{'})
    \\
    0   & b_{32}^{-}(\xi^{'}) &  0 \ep,
  \ea 
 Combining Lemma \ref{Upper bounds for the symbolic flow $2*2$} and Remark \ref{Remark 2.1}, we can obtain  \eqref{(2.14)} for  $\xi\in\mathcal{R}_{12}\backslash(\mathcal{R}_{23}+k)$ and $\xi\in(\mathcal{R}_{23}+k)\backslash\mathcal{R}_{12}$ directly.
 
  While in $\mathcal{R}_{12}\cap(\mathcal{R}_{23}+k)$, two non-transparent resonances $(1,2)$ and $(2,3)$ are coupled together. 
  By Assumption \ref{Relaxed  separation conditions}, if at least one of $b_{32}^{-}$, $b_{23}^{+}$ is transparent in $(\mathcal{R}_{12}-k)\cap\mathcal{R}_{23}$, then at least one of $(b_{32}^{-})_{-1}$, $(b_{23}^{+})_{-1}$ is transparent in $\mathcal{R}_{12}\cap(\mathcal{R}_{23}+k)$. That is to say one of $(b_{32}^{-})_{-1}$, $(b_{23}^{+})_{-1}$ is transparent in $\mathcal{R}_{12}\cap(\mathcal{R}_{23}+k)$ or both  $(b_{32}^{-})_{-1}$, $(b_{23}^{+})_{-1}$ are transparent in $\mathcal{R}_{12}\cap(\mathcal{R}_{23}+k)$, then we can eliminate transparent term in $\mathcal{R}_{12}\cap(\mathcal{R}_{23}+k)$ and similarly non-transparent properties of interaction coefficients $b_{12}^{+}$, $b_{21}^{-}$ in $\mathcal{R}_{12}\cap(\mathcal{R}_{23}+k)$ are considered. Under different separation conditions,  \eqref{(2.13)} can be divided into different equations.

  Firstly, if at least one of $b_{32}^{-}$, $b_{23}^{+}$ is transparent in $(\mathcal{R}_{12}-k)\cap\mathcal{R}_{23}$ and there is only $b_{23}^{+}$ transparent in $(\mathcal{R}_{12}-k)\cap\mathcal{R}_{23}$. Using normal form reduction, we have
    \ba
  \check{B}=\chi_{12}\chi_{23,-1}
  \bp  0 & b_{12}^{+}  & 0
    \\ 
     b_{21}^{-} & 0 &0
    \\
    0   & (b_{32}^{-})_{-1} &  0 \ep.\nn
  \ea
Let 
\be\label{(2.20)}
V=\op_{\e}(\chi_{12}\chi_{23,-1})(\varphi_{12}\varphi_{23}  \check{U}),
\ee
\be\label{(2.21)}
W=(W_{1},W_{2})=(\op_{\e}(\chi_{12}\chi_{23,-1})((1-\varphi_{12}\varphi_{23} ) \check{U}), (1-\op_{\e}(\chi_{12}\chi_{23,-1}))\check{U})
\ee
satisfying,
\be
\check{U}=V+W_{1}+W_{2}.\nn
\ee
Then we obtain the following equations:
\be\label{(2.22)}
\partial_{t}V+\frac{1}{\sqrt{\e}}\op_{\e}^{\psi}(M)V=\sqrt{\e}F_{V},
\ee
\be\label{(2.23)}
\partial_{t}W+\frac{1}{\sqrt{\e}}\op_{\e}(iA)W=\op_{\e}(D)W+\sqrt{\e}F_{W}
\ee
with symbol
\be\label{(2.24)}
M=i(\chi_{12}\chi_{23,-1})^{\#} A-\sqrt{\e}(\varphi_{12}\varphi_{23})^{\#}  \check{B},\quad
 A=\bp \mathcal{A} & 0
 \\
 0 &\mathcal{A}\ep,\quad
 D=(1-(\varphi_{12}\varphi_{23})^{\#})(\chi_{12}\chi_{23,-1})^{\#}\check{B}.
 \ee
Based on the given $M$, we now consider the equation
\be\label{(2.25)}
\partial_{t}S_{0}+\frac{1}{\sqrt{\e}}MS_{0}=0, \quad S_{0}(\tau,\tau)=\Id.
\ee
Then combine \eqref{(2.24)} to obtain
  \ba\label{M-form}
  M&=\bp 
  i(\chi_{12}\chi_{23,-1})^{\#}(\l_{1,+1}-\omega) & -\sqrt{\e}(\varphi_{12}\varphi_{23})^{\#}\chi_{12}\chi_{23,-1}b_{12}^{+}  & 0
    \\ 
     -\sqrt{\e}(\varphi_{12}\varphi_{23})^{\#}\chi_{12}\chi_{23,-1}b_{21}^{-} & i(\chi_{12}\chi_{23,-1})^{\#}\l_{2} &0    \\
    0   &-\sqrt{\e}(\varphi_{12}\varphi_{23})^{\#}\chi_{12}\chi_{23,-1}(b_{32}^{-})_{-1} &  i(\chi_{12}\chi_{23,-1})^{\#}(\l_{3,-1}+\omega) \ep\\
        &=:\bp 
  i\mu_{1}^{'} & -\sqrt{\e}(b_{12}^{+})^{'} & 0
    \\ 
     -\sqrt{\e}(b_{21}^{-})^{'} & i\mu_{2}^{'} &0  
      \\
    0   &-\sqrt{\e}(b_{32}^{-})_{-1}^{'} &  i\mu_{3}^{'} \ep. 
  \ea
 We denote $S_{0}(\tau,t)=\bp y_{1}\\ y_{2}\\y_{3}\ep(t)$, then \eqref{(2.25)} can be write into
  \ba\label{(2.26)}
\partial_{t}\bp y_{1}\\ y_{2}\ep+\frac{1}{\sqrt{\e}}
M_{12}\bp y_{1}\\ y_{2}\ep=0, 
\ea
\be\label{(2.27)}
\partial_{t}y_{3}-(b_{32}^{-})_{-1}^{'}y_{2}+\frac{1}{\sqrt{\e}}i\mu_{3}^{'}y_{3}=0,
\ee
with
    \ba
    M_{12}=\bp 
  i\mu_{1}^{'} & -\sqrt{\e}(b_{12}^{+})^{'}
    \\ 
     -\sqrt{\e}(b_{21}^{-})^{'} & i\mu_{2}^{'}\ep. \nn
  \ea  
 For the case of \eqref{(2.26)}-\eqref{(2.27)}, by Lemma \ref{Upper bounds for the symbolic flow $2*2$}, we have
\ba
|\bp y_{1}\\ y_{2}\ep(t)|=|e^{\frac{\tau-t}{\sqrt{\e}}M_{12}}\cdot \bp y_{1}\\ y_{2}\ep(\tau)|\lesssim\exp((t-\tau){\g}_{12}^{+}). \nn
\ea
Then for $|y_{3}|$, we have
\ba
|y_{3}(t)|&=|e^{\frac{(\tau -t)}{\sqrt{\e}}i\mu_{3}^{'}}y_{3}(\tau)+\int_{\tau}^{t}e^{(t^{'}-t)\cdot\frac{1}{\sqrt{\e}}i\mu_{3}}\cdot (b_{32}^{-})_{-1}^{'}y_{2}(t^{'})dt^{'}|\nn
\\
&\lesssim|e^{\frac{(\tau -t)}{\sqrt{\e}}i\mu_{3}^{'}}|\cdot|y_{3}(\tau)|+\int_{\tau}^{t}|e^{(t^{'}-t)\cdot\frac{1}{\sqrt{\e}}i\mu_{3}}|\cdot |(b_{32}^{-})_{-1}^{'}|\cdot|y_{2}(t^{'})|dt^{'}\\
&\lesssim 1+CT|\ln\e|\exp((t-\tau){\g}_{12}^{+})\lesssim \exp(t{\g}_{12}^{+}).
\ea
Finally, we obtain
\be\label{(2.28)}
|S_{0}(\tau,t)|=|\bp y_{1}\\ y_{2}\\y_{3}\ep(t)|\lesssim \exp(t\g_{12}^{+})\lesssim\exp(t\g^{+}).
\ee
Secondly, if at least one of $b_{32}^{-}$, $b_{23}^{+}$ is transparent in $(\mathcal{R}_{12}-k)\cap\mathcal{R}_{23}$ and there is only $b_{32}^{-}$ transparent in $(\mathcal{R}_{12}-k)\cap\mathcal{R}_{23}$. After normal form reduction, we have
  \ba
  \check{B}=\chi_{12}\chi_{23,-1}
  \bp  0 & b_{12}^{+}  & 0
    \\ 
     b_{21}^{-} & 0 &(b_{23}^{+})_{-1}
    \\
    0   & 0 &  0 \ep.\nn
  \ea
  
 and
  \ba
  M&=\bp 
  i(\chi_{12}\chi_{23,-1})^{\#}(\l_{1,+1}-\omega) & -\sqrt{\e}(\varphi_{12}\varphi_{23})^{\#}\chi_{12}\chi_{23,-1}b_{12}^{+}  & 0
    \\ 
     -\sqrt{\e}(\varphi_{12}\varphi_{23})^{\#}\chi_{12}\chi_{23,-1}b_{21}^{-} & i(\chi_{12}\chi_{23,-1})^{\#}\l_{2} &-\sqrt{\e}(\varphi_{12}\varphi_{23})^{\#}\chi_{12}\chi_{23,-1}(b_{23}^{+})_{-1}
    \\
    0   &0 &  i(\chi_{12}\chi_{23,-1})^{\#}(\l_{3,-1}+\omega) \ep\\
        &=:\bp 
  i\mu_{1}^{'} & -\sqrt{\e}(b_{12}^{+})^{'} & 0
    \\ 
     -\sqrt{\e}(b_{21}^{-})^{'} & i\mu_{2}^{'} &-\sqrt{\e}(b_{23}^{+})_{-1}^{'}  
      \\
    0   &0 &  i\mu_{3}^{'} \ep. \nn
  \ea
   Then \eqref{(2.25)} can be write into
\ba\label{(2.29)}
\partial_{t}\bp y_{1}\\ y_{2}\ep+\frac{1}{\sqrt{\e}}M_{12}\bp y_{1}\\ y_{2}\ep=\bp 0\\ (b_{23}^{+})_{-1}^{'}y_{3}\ep,
\ea
\be\label{(2.30)}
\partial_{t}y_{3}+\frac{1}{\sqrt{\e}}i\mu_{3}^{'}y_{3}=0.
\ee
For \eqref{(2.29)}-\eqref{(2.30)}, we have
\be
\bp y_{1}\\ y_{2}\ep(t)=e^{\frac{\tau-t}{\sqrt{\e}}M_{12}}\bp y_{1}\\ y_{2}\ep(\tau)+\int_{\tau}^{t}e^{(t^{'}-t)\cdot\frac{1}{\sqrt{\e}}M_{12}}\cdot 
\bp 0\\ (b_{23}^{+})_{-1}^{'}y_{3}\ep(t^{'})dt^{'}, \nn
\ee
\be
y_{3}(t)=e^{\frac{\tau-t}{\sqrt{\e}}i\mu_{3}^{'}}.\nn
\ee
It is obvious that
\be
|y_{3}|\lesssim 1\nn
\ee
and by Lemma \ref{Upper bounds for the symbolic flow $2*2$}, we have
\ba
|\bp y_{1}\\ y_{2}\ep(t)|&\lesssim |e^{\frac{\tau-t}{\sqrt{\e}}M_{12}}|\cdot|\bp y_{1}\\ y_{2}\ep(\tau)|+\int_{\tau}^{t}|e^{(t^{'}-t)\cdot\frac{1}{\sqrt{\e}}M_{12}}|\cdot |\bp 0\\ (b_{23}^{+})_{-1}^{'}y_{3}\ep(t^{'})|dt^{'}\\
&\lesssim \exp((t-\tau){\g}_{12}^{+})\cdot 1+CT|\ln\e|\exp(t{\g}_{12}^{+})\lesssim \exp(t{\g}_{12}^{+}).\nn
\ea
Thus we have,
\be\label{(2.31)}
|S_{0}(\tau,t)|=|\bp y_{1}\\ y_{2}\\y_{3}\ep|\lesssim \exp(t{\g_{12}}^{+})\lesssim\exp(t{\g}^{+}).
\ee

 Similarly, if at least one of $b_{21}^{-}$, $b_{12}^{+}$ is transparent in $\mathcal{R}_{12}\cap(\mathcal{R}_{23}+k)$, we can also obtain the same upper bounds for $S_{0}$. Finally, if both  $(b_{32}^{-})_{-1}$, $(b_{23}^{+})_{-1}$ are transparent in $\mathcal{R}_{12}\cap(\mathcal{R}_{23}+k)$, after normal form reduction,  $M$ in \eqref{(2.25)} can be reduced to $M_{12}$ or if both $b_{12}^{+}$, $b_{21}^{-}$ are transparent in $\mathcal{R}_{12}\cap(\mathcal{R}_{23}+k)$, $M$ in \eqref{(2.25)} can be reduced to $M_{23}$. Thus we can solve \eqref{(2.25)} respectively with $|S_{0}(\tau,t)|\lesssim\exp(t\g_{12}^{+})$ and $|S_{0}(\tau,t)|\lesssim\exp(t\g_{23}^{+})$ directly using Lemma \ref{Upper bounds for the symbolic flow $2*2$}.
\subsection{Case 2.}\label{Case 2.}
Assume $(1,2)$, $(1,3)\in \Re_{0}$, at least one of  $b_{13}^{+}$, $b_{31}^{-}$ is transparent on $\mathcal {R}_{12}\bigcap\mathcal {R}_{13}$, or at least one of  $b_{12}^{+}$, $b_{21}^{-}$ is transparent on $\mathcal {R}_{12}\bigcap\mathcal {R}_{13}$.
 By the eigenmodes of the hyperbolic operator, we decompose $\dot u$ and shift the component related to $\Pi_{1}$. Then we define 
\be
U_{1}=e^{-i\th}\op_{\e}{(\Pi_1)}\dot u, \quad U_{2}=\op_{\e}{(\Pi_2)}\dot u, \quad U_{3}=\op_{\e}{(\Pi_3)}\dot u, \quad  U_{4}=\op_{\e}{(\Pi_4)}\dot u,\quad\cdots,\quad U_{J}=\op_{\e}{(\Pi_J)}\dot u.\nn
\ee
The perturbation unknown $\dot u$ is decomposed by
\be
\dot u= e^{i\th}U_{1}+U_{2}+U_{3}+U_{4}+\cdots+U_{J},\nn
\ee
and
\be
A=\diag(\l_{1,+1}-w,\l_{2},\l_{3},\cdots,\l_{J}).\nn
\ee
 Similar to Case 1, we utilize normal form reduction to have
\ba
  \check{B}=\bp  0 &  \chi_{12}b_{12}^{+}  & \chi_{13}b_{13}^{+}
    \\ 
     \chi_{12}b_{21}^{-} & 0 & 0
     \\
     \chi_{13}b_{31}^{-}  & 0 &  0 \ep,\nn
  \ea
and when $\xi\in \mathcal{R}_{12}\cap\mathcal{R}_{13}$ for $\mathcal{R}_{12}\cap\mathcal{R}_{13}\neq \varnothing$, we have
\ba\label{(2.32)}
M&=i(\chi_{12}\chi_{13})^{\#} A-\sqrt{\e}(\varphi_{12}\varphi_{13})^{\#}  \check{B}\\
   &=\bp 
  i(\chi_{12}\chi_{13})^{\#}(\l_{1,+1}-w) & -\sqrt{\e}(\varphi_{12}\varphi_{13})^{\#}\chi_{12}\chi_{13}b_{12}^{+}  & -\sqrt{\e}(\varphi_{12}\varphi_{13})^{\#}\chi_{12}\chi_{13}b_{13}^{+}
    \\ 
     -\sqrt{\e}(\varphi_{12}\varphi_{13})^{\#}\chi_{12}\chi_{13}b_{21}^{-} & i(\chi_{12}\chi_{13})^{\#}\l_{2} &0
    \\
    -\sqrt{\e}(\varphi_{12}\varphi_{13})^{\#}\chi_{12}\chi_{13}b_{31}^{-}&0 &  i(\chi_{12}\chi_{13})^{\#}\l_{3} \ep.
 \ea
  Here the method of variable substitution and the classified discussion on intersections of resonance sets occupy an essential position. 
  
  We now consider the solutions of $\partial_{t}S_{0}+\frac{1}{\sqrt{\e}}MS_{0}=0, S_{0}(\tau,\tau)=\Id$ in different non-transparent resonance sets $\mathcal{R}_{12}$ and $\mathcal{R}_{13}$.  When $\xi\in \mathcal{R}_{12}$ for $\mathcal{R}_{12}\cap\mathcal{R}_{13}= \varnothing$, $\xi\in \mathcal{R}_{12}\backslash \mathcal{R}_{13}$ for $\mathcal{R}_{12}\cap\mathcal{R}_{13}\neq \varnothing$ and $\xi\in \mathcal{R}_{12}\cap\mathcal{R}_{13}$ for $\mathcal{R}_{12}\cap\mathcal{R}_{13}\neq \varnothing$ with the condition that at least one of  $b_{13}^{+}$, $b_{31}^{-}$ is transparent in $\mathcal {R}_{12}\bigcap\mathcal {R}_{13}$ and utilizing normal form reduction to eliminate the transparent terms, we can obtain the following bound for the symbolic flow 
 \be\label{(2.33)}
 |S_{0}(\tau,t)|\lesssim|\ln\e|^{*}\exp(t\gamma_{12}^{+})\lesssim|\ln\e|^{*}\exp(t\gamma^{+}).
 \ee
  Moreover, when $\xi\in \mathcal{R}_{13}$ for $\mathcal{R}_{12}\cap\mathcal{R}_{13}= \varnothing$, $\xi\in \mathcal{R}_{13}\backslash \mathcal{R}_{12}$ for $\mathcal{R}_{12}\cap\mathcal{R}_{13}\neq \varnothing$ and $\xi\in \mathcal{R}_{12}\cap\mathcal{R}_{13}$ for $\mathcal{R}_{12}\cap\mathcal{R}_{13}\neq \varnothing$ with the condition that at least one of  $b_{12}^{+}$, $b_{21}^{-}$ is transparent on $\mathcal {R}_{12}\bigcap\mathcal {R}_{13}$, we similarly have
  \be\label{(2.34)}
  |S_{0}(\tau,t)|\lesssim|\ln\e|^{*}\exp(t\gamma_{13}^{+})\lesssim|\ln\e|^{*}\exp(t\gamma^{+}).
  \ee
\subsection{Case 3}\label{Case 3.}
\quad In Case 3, for resonances $(1,3), (2,3)\in \Re_{0}$, at least one of $b_{23}^{+}$, $b_{32}^{-}$ is transparent on $\mathcal {R}_{13}\bigcap\mathcal {R}_{23}$, or at least one of $b_{13}^{+}$, $b_{31}^{-}$ is transparent on $\mathcal {R}_{13}\bigcap\mathcal {R}_{23}$, we decompose $\dot u$ and shift the component related to  $\Pi_{3}$ and define 
  \be
U_{1}=\op_{\e}{(\Pi_1)}\dot u, \quad U_{2}=\op_{\e}{(\Pi_2)}\dot u,\quad  U_{3}=e^{i\th}\op_{\e}{(\Pi_3)}\dot u, \quad U_{4}=\op_{\e}{(\Pi_4)}\dot u,\quad \cdots, \quad U_{J}=\op_{\e}{(\Pi_J)}\dot u,\nn
\ee
and
\be
A=\diag(\l_1,\l_{2},\l_{3,-1}+w,\l_{4},\cdots,\l_{J}).\nn
\ee
Then, similar to Case 1 and Case 2, we have
\ba\label{(2.35)}
  \check{B}=\bp  0 &  0  & \chi_{13,-1}(b_{13}^{+})_{-1}
    \\ 
     0  & 0 & \chi_{23,-1}(b_{23}^{+})_{-1}
     \\
     \chi_{13,-1}(b_{31}^{-})_{-1}  &  \chi_{23,-1}(b_{32}^{-})_{-1} &  0 \ep(\xi),
  \ea
  and when $\xi\in (\mathcal{R}_{13}+k)\cap (\mathcal{R}_{23}+k)$ for $(\mathcal{R}_{13}+k)\cap(\mathcal{R}_{23}+k)\neq \varnothing$, we have
\ba
  M&=i(\chi_{13,-1}\chi_{23,-1})^{\#} A-\sqrt{\e}(\varphi_{13}\varphi_{23})^{\#} \check{B}\\
&=\small{i(\chi_{13,-1}\chi_{23,-1})^{\#}
\bp  \l_{1} &  0  & 0
    \\ 
     0 & \l_{2} & 0
    \\
    0   & 0 &  \l_{3,-1}+\omega \ep
    -\sqrt{\e}(\varphi_{13}\varphi_{23})^{\#} (\chi_{13,-1}\chi_{23,-1})
    \bp  0 &  0  &(b_{13}^{+})_{-1}
    \\ 
     0  & 0 & (b_{23}^{+})_{-1}
     \\
     (b_{31}^{-})_{-1}  & (b_{32}^{-})_{-1} &  0 \ep}\\
    &=\small{ \bp  i(\chi_{13,-1}\chi_{23,-1})^{\#}\l_{1} &  0  & \sqrt{\e}(\varphi_{13}\varphi_{23})^{\#} (\chi_{13,-1}\chi_{23,-1})(b_{13}^{+})_{-1}
    \\ 
     0 &i(\chi_{13,-1}\chi_{23,-1})^{\#}\l_{2} & \sqrt{\e}(\varphi_{13}\varphi_{23})^{\#} (\chi_{13,-1}\chi_{23,-1})(b_{23}^{+})_{-1}
    \\
    \sqrt{\e}(\varphi_{13}\varphi_{23})^{\#} (\chi_{13,-1}\chi_{23,-1})(b_{31}^{-})_{-1}   & \sqrt{\e}(\varphi_{13}\varphi_{23})^{\#} (\chi_{13,-1}\chi_{23,-1}) (b_{32}^{-})_{-1}&  i(\chi_{13,-1}\chi_{23,-1})^{\#}(\l_{3,-1}+\omega) \ep}\\
    &=:\small{ \bp 
  i\mu_{1}^{'} & 0 & -\sqrt{\e}(b_{13}^{+})_{-1}^{'}
    \\ 
     0 & i\mu_{2}^{'} &-\sqrt{\e}(b_{23}^{+})_{-1}^{'}  
      \\
    -\sqrt{\e}(b_{31}^{-})_{-1}^{'}   &-\sqrt{\e}(b_{32}^{-})_{-1}^{'} &  i\mu_{3}^{'} \ep}.\nn
\ea

Combining Remark \ref{Remark 2.1},  when $\xi\in (\mathcal{R}_{13}+k)$ for $(\mathcal{R}_{13}+k)\cap(\mathcal{R}_{23}+k)= \varnothing$, $\xi\in (\mathcal{R}_{13}+k)\backslash (\mathcal{R}_{23}+k)$ for $(\mathcal{R}_{13}+k)\cap(\mathcal{R}_{23}+k)\neq\varnothing$ and 
$\xi\in (\mathcal{R}_{13}+k)\cap (\mathcal{R}_{23}+k)$ for $(\mathcal{R}_{13}+k)\cap(\mathcal{R}_{23}+k)\neq \varnothing$  with the condition that at least one of  $b_{23}^{+}$, $b_{32}^{-}$ is transparent in $\mathcal {R}_{12}\bigcap\mathcal {R}_{13}$, we can obtain following bound for the  symbolic flow 
 \be\label{(2.36)}
 |S_{0}(\tau,t)|\lesssim|\ln\e|^{*}\exp(t\gamma_{13}^{+})\lesssim|\ln\e|^{*}\exp(t\gamma^{+}).
 \ee
  Moreover, when $\xi\in (\mathcal{R}_{23}+k)$ for $(\mathcal{R}_{13}+k)\cap(\mathcal{R}_{23}+k)= \varnothing$, $\xi\in (\mathcal{R}_{23}+k)\backslash (\mathcal{R}_{13}+k)$ for $(\mathcal{R}_{13}+k)\cap(\mathcal{R}_{23}+k)\neq\varnothing$ and $\xi\in (\mathcal{R}_{13}+k)\cap (\mathcal{R}_{23}+k)$ for $(\mathcal{R}_{13}+k)\cap(\mathcal{R}_{23}+k)\neq \varnothing$  with the condition that at least one of  $b_{13}^{+}$, $b_{31}^{-}$ is transparent in $\mathcal {R}_{13}\bigcap\mathcal {R}_{23}$, we similarly have
  \be\label{(2.37)}
  |S_{0}(\tau,t)|\lesssim|\ln\e|^{*}\exp(t\gamma_{23}^{+})\lesssim|\ln\e|^{*}\exp(t\gamma^{+}).
  \ee
  
 \subsection{Instability of WKB solution} \label{Instability of WKB solution}
Combining Case 1 to Case 3, for all $T>0$, $0\leq \tau\leq t\leq T|\ln\e|$, and under the given Assumptions in Theorem \ref{New instability criterion}, we can conclude that the solution $S_{0}$ to 
 \be
\partial_{t}S_{0}+\frac{1}{\sqrt{\e}}M_{ijj^{'}}S_{0}=0, \quad S_{0}(\tau,\tau)=\Id\nn
\ee
always satisfies the upper bounds
\be
|S_{0}(\tau,t)|\lesssim|\ln\e|^{*}\exp(t\gamma^{+}),\nn
\ee
where $\{M_{ijj^{'}}\}_{1\leq i,j,j^{'}\leq J}$ are $3\times 3$ block matrices of non-transparent resonance terms after normal form reductions.
  Similar to analysis of general cases in Section 6.3.2 of \cite{Lu}, we have the following estimate:
\be
|\partial_{x}^{\alpha}S_{0}(\tau,t)|\lesssim|\ln\e|^{*}\exp(t\gamma^{+}).\nn
\ee
 Combining \cite{Lu}, we prove the first result in Theorem \ref{New instability criterion} using the same method of Duhamel representation for the instability. 

 \section{Three coupled resonances}\label{Three coupled resonances}
 For sake of simplicity, it takes no influences if we directly consider non-transparent resonances set $\Re_{1}=\{(1,2), (2,3), (1,3)\}\subset \Re_{0}$. We have the freedom to choose different frequency shift and we will see that we can do some flexible extensions for the separation conditions in the arguments. By the eigenmodes of the hyperbolic operator, we decompose $\dot u$ and shift the component related to $\Pi_{1}$ as follows 
\be
U_{1}=e^{-i\th}\op_{\e}{(\Pi_1)}\dot u,\quad U_{2}=\op_{\e}{(\Pi_2)}\dot u,\quad U_{3}=\op_{\e}{(\Pi_3)}\dot u,\quad U_{4}=\op_{\e}{(\Pi_4)}\dot u,\quad\cdots, \quad U_{J}=\op_{\e}{(\Pi_J)}\dot u.\nn
\ee
The perturbation unknown $\dot u$ is decomposed by
\be
\dot u= e^{i\th}U_{1}+U_{2}+U_{3}+\cdots+U_{J}\nn
\ee
and
\be
A=\diag(\l_{1,+1}-w,\l_{2},\l_{3},\cdots,\l_{J}).\nn
\ee

Then, similar to Case 1 to Case 3 in Section \ref{Cases of two coupled resonances}, we have
\ba
  B=\sum_{p=\pm 1}
  \bp  e^{ip\th}\Pi_{1,+(p+1)}B_{p}\Pi_{1,+1} &  e^{i(p-1)\th}\Pi_{1,+p}B_{p}\Pi_{2}  &  e^{i(p-1)\th}\Pi_{1,+p}B_{p}\Pi_{3}&\cdots& e^{i(p-1)\th}\Pi_{1,+p}B_{p}\Pi_{J}
    \\ 
     e^{i(p+1)\th}\Pi_{2,+(p+1)}B_{p}\Pi_{1,+1} & e^{ip\th}\Pi_{2,+p}B_{p}\Pi_{2} & e^{ip\th}\Pi_{2,+p}B_{p}\Pi_{3}&\cdots& e^{ip\th}\Pi_{2,+p}B_{p}\Pi_{J}
    \\
     e^{i(p+1)\th}\Pi_{3,+(p+1)}B_{p}\Pi_{1,+1} & e^{ip\th}\Pi_{3,+p}B_{p}\Pi_{2} &  e^{ip\th}\Pi_{3,+p}B_{p}\Pi_{3} &\cdots&e^{ip\th}\Pi_{3,+p}B_{p}\Pi_{J}
 \\
 \vdots &\vdots &\vdots & &\vdots 
 \\
     e^{i(p+1)\th}\Pi_{J,+(p+1)}B_{p}\Pi_{1,+1} & e^{ip\th}\Pi_{J,+p}B_{p}\Pi_{2} &  e^{ip\th}\Pi_{J,+p}B_{p}\Pi_{3} &\cdots&e^{ip\th}\Pi_{J,+p}B_{p}\Pi_{J}\ep.\nn
  \ea
Since we only have $(1,2), (1,3), (2,3)\in \Re_{0}$, which means $(3,2)\notin \Re_{0}$. Then after normal form reduction,  we have the following $\check B$ with oscillation terms:
\ba\label{(2.38)}
  \check{B}=\bp  0 &  \chi_{12}b_{12}^{+}  & \chi_{13}b_{13}^{+}
    \\ 
     \chi_{12}b_{21}^{-} & 0 & e^{i\th}\chi_{23}b_{23}^{+}
     \\
     \chi_{13}b_{31}^{-}  & e^{-i\th}\chi_{23,-1}(b_{32}^{-})_{-1} &  0 \ep.
  \ea

Since normal form reduction can help us to eliminate transparent interaction coefficients in above $B$, here we assume $\mathcal {R}_{12}\bigcap\mathcal {R}_{13}\bigcap\mathcal {R}_{23}\neq \varnothing$ and $b_{23}^{+}$ is transparent in $\mathcal {R}_{23}$, then $\mathcal {R}_{12}\bigcap\mathcal {R}_{13}\bigcap(\mathcal {R}_{23}+k)= \varnothing$. Or, we assume $\mathcal {R}_{12}\bigcap\mathcal {R}_{13}\bigcap(\mathcal {R}_{23}+k)\neq \varnothing$ and $b_{32}^{-}$ is transparent in $\mathcal {R}_{23}$, then $\mathcal {R}_{12}\bigcap\mathcal {R}_{13}\bigcap\mathcal {R}_{23}= \varnothing$. Thus we can eliminate the oscillation terms $e^{i\th}\chi_{23}b_{23}^{+}$ and $e^{-i\th}\chi_{23,-1}(b_{32}^{-})_{-1}$ in $\check{B}$. And combining the assumption that one of  $b_{13}^{+}$, $b_{31}^{-}$ is transparent in $\mathcal {R}_{13}$, or one of  $b_{12}^{+}$, $b_{21}^{-}$ is transparent in $\mathcal {R}_{12}$, we eliminate one of   $\chi_{13}b_{13}^{+}$, $\chi_{13}b_{31}^{-}$, $\chi_{12}b_{12}^{+}$, $\chi_{12}b_{21}^{-}$ in $\check{B}$ in normal form reduction and we can obtain the  coupled case of $\check{B}$ as follows 
\ba
  \check{B}_{1}=\bp  0 & 0 & \chi_{13}b_{13}^{+}
    \\ 
     \chi_{12}b_{21}^{-} & 0 & 0
     \\
     \chi_{13}b_{31}^{-}  &0 &  0 \ep,
       \check{B}_{2}=\bp  0 &  \chi_{12}b_{12}^{+}  & 0
    \\ 
     \chi_{12}b_{21}^{-} & 0 & 0
     \\
     \chi_{13}b_{31}^{-}  &0 &  0 \ep,\nn
  \ea
  \ba
  \check{B}_{3}=\bp  0 &  \chi_{12}b_{12}^{+}  & \chi_{13}b_{13}^{+}
    \\ 
0 & 0 & 0
     \\
     \chi_{13}b_{31}^{-}  &0 &  0 \ep,
       \check{B}_{4}=\bp  0 &  \chi_{12}b_{12}^{+}  & \chi_{13}b_{13}^{+}
    \\ 
     \chi_{12}b_{21}^{-} & 0 & 0
     \\
    0  &0 &  0 \ep.\nn
  \ea
Now we use the combine Case 2 in Section \ref{Cases of two coupled resonances}, we obtain the bound of symmetric flow: $|S_{0}(\tau,t)|\lesssim \exp(t{\g}^{+})$. 

Following the way to Case 1 in Section \ref{Cases of two coupled resonances} and by the eigenmodes of hyperbolic operator, we decompose $\dot u$ and shift the component related to $\Pi_{1}$ and  $\Pi_{3}$. Then we can make the following frequency shift 
\be
U_{1}=e^{-i\th}\op_{\e}{(\Pi_1)}\dot u,\quad  U_{2}=\op_{\e}{(\Pi_2)}\dot u,\quad  U_{3}=e^{i\th}\op_{\e}{(\Pi_3)}\dot u,U_{4}=\op_{\e}{(\Pi_4)}\dot u,\quad\cdots, \quad U_{J}=\op_{\e}{(\Pi_J)}\dot u.\nn
\ee
 Then the perturbation unknown $\dot u$ is decomposed by
\be
\dot u=e^{i\th}U_{1}+U_{2}+e^{-i\th}U_{3}+U_{4}+\cdots+U_{J},\nn
\ee
and the symbol of the propagator is 
\be
A=\diag(\l_{1,+1}-w,\l_{2},\l_{3,-1}+\omega,\cdots,\l_{J}).\nn
\ee

The symbol of singular source term is
  \ba
  B=\small{\sum_{p=\pm 1}
  \bp  e^{ip\th}\Pi_{1,+(p+1)}B_{p}\Pi_{1,+1} &  e^{i(p-1)\th}\Pi_{1,+p}B_{p}\Pi_{2}  &  e^{i(p-2)\th}\Pi_{1,+(p-1)}B_{p}\Pi_{3,-1}&\cdots &  e^{i(p-1)\th}\Pi_{1,+p}B_{p}\Pi_{J}
    \\ 
     e^{i(p+1)\th}\Pi_{2,+(p+1)}B_{p}\Pi_{1,+1} & e^{ip\th}\Pi_{2,+p}B_{p}\Pi_{2} & e^{i(p-1)\th}\Pi_{2,+(p-1)}B_{p}\Pi_{3,-1}&\cdots &  e^{ip\th}\Pi_{2,+p}B_{p}\Pi_{J}
    \\
     e^{i(p+2)\th}\Pi_{3,+(p+1)}B_{p}\Pi_{1,+1} & e^{i(p+1)\th}\Pi_{3,+p}B_{p}\Pi_{2} &  e^{ip\th}\Pi_{3,+(p-1)}B_{p}\Pi_{3,-1}&\cdots &  e^{i(p+1)\th}\Pi_{3,+p}B_{p}\Pi_{J}
     \\ 
     \vdots & \vdots & \vdots& & \vdots
     \\
     e^{i(p+1)\th}\Pi_{J,+(p+1)}B_{p}\Pi_{1,+1} &  e^{ip\th}\Pi_{J,+(p-1)}B_{p}\Pi_{2,-1}  &  e^{i(p-1)\th}\Pi_{J,+(p-1)}B_{p}\Pi_{3,-1}&\cdots &  e^{ip\th}\Pi_{J,+p}B_{p}\Pi_{J}\ep}\nn
     \ea
  with $\Pi_{j,+q}{(\xi)}=\Pi(\xi+qk)$ for $q\in\ZZ$ and 
  \be
  B_{p}=B(u_{0,p}), \quad p\in{\{-1,1\}},\nn
  \ee
  where $u_{0,\pm 1}$ are the leading amplitudes in the WKB solution. Combining the separation conditions in Assumption \ref{Relaxed  separation conditions}, after normal form reduction we have
\ba\label{(2.39)}
  \check{B}&=\bp  0 &  \chi_{12}\Pi_{1,+1}B_{1}\Pi_{2}  & \chi_{13,-1}e^{-i\th}\Pi_{1}B_{1}\Pi_{3,-1}
    \\ 
     \chi_{12}\Pi_{2}B_{-1}\Pi_{1,+1} & 0 & \chi_{23,-1}\Pi_{2}B_{1}\Pi_{3,-1}
    \\
    \chi_{13}e^{i\th}\Pi_{3}B_{-1}\Pi_{1,+1}   & \chi_{23,-1}\Pi_{3,-1}B_{-1}\Pi_{2} &  0 \ep\\
   & =\bp  0 &  \chi_{12}b_{12}^{+}  & \chi_{13,-1}e^{-i\th}(b_{13}^{+})_{-1}
    \\ 
     \chi_{12}b_{21}^{-} & 0 & \chi_{23,-1}(b_{23}^{+})_{-1}
    \\
    \chi_{13}e^{i\th}b_{31}^{-}    & \chi_{23,-1}(b_{32}^{-})_{-1} &  0 \ep.
  \ea  
  
Assume $\mathcal {R}_{12}\bigcap(\mathcal {R}_{13}+k)\bigcap(\mathcal {R}_{23}+k)\neq \varnothing$ and $b_{13}^{+}$ is transparent in $\mathcal {R}_{13}$, then $\mathcal {R}_{12}\bigcap\mathcal {R}_{13}\bigcap(\mathcal {R}_{23}+k)=\varnothing$ shows the oscillation terms $\chi_{13,-1}e^{-i\th}(b_{13}^{+})_{-1}$, $\chi_{13}e^{i\th}b_{31}^{-}$ will disappear in normal form reduction. Or, we assume $\mathcal {R}_{12}\bigcap\mathcal {R}_{13}\bigcap(\mathcal {R}_{23}+k)\neq\varnothing$ and $b_{31}^{-}$ is transparent in $\mathcal {R}_{13}$, then $\mathcal {R}_{12}\bigcap(\mathcal {R}_{13}+k)\bigcap(\mathcal {R}_{23}+k)= \varnothing$ will also eliminate the oscillation terms in normal form reduction. Combining the assumption that one of  $b_{12}^{+}$, $b_{21}^{-}$ is transparent in $\mathcal {R}_{13}$, or one of  $b_{23}^{+}$, $b_{32}^{-}$ is transparent in $\mathcal {R}_{23}$, the coupled case of $\check{B}$ could be reduced to one of the following cases 
\ba
  \check{B}=\bp  0 &  0  & 0
    \\ 
     \chi_{12}b_{21}^{-} & 0 & \chi_{23,-1}(b_{23}^{+})_{-1}
    \\
    0& \chi_{23,-1}(b_{32}^{-})_{-1} &  0 \ep,
  \check{B}=\bp  0 &  \chi_{12}b_{12}^{+}  & 0
    \\ 
     0 & 0 & \chi_{23,-1}(b_{23}^{+})_{-1}
    \\
    0& \chi_{23,-1}(b_{32}^{-})_{-1} &  0 \ep,\nn
  \ea
  \ba
  \check{B}=\bp  0 &  \chi_{12}b_{12}^{+}  & 0
    \\ 
     \chi_{12}b_{21}^{-} & 0 & 0
    \\
    0& \chi_{23,-1}(b_{32}^{-})_{-1} &  0 \ep,
  \check{B}=\bp  0 &  \chi_{12}b_{12}^{+}  & 0
    \\ 
     \chi_{12}b_{21}^{-} & 0 & \chi_{23,-1}(b_{23}^{+})_{-1}
    \\
    0& 0 &  0 \ep.\nn
  \ea
Combine the cases 1  in Section \ref{Cases of two coupled resonances}, we can obtain the bound of symmetric flows: $|S_{0}(\tau,t)|\lesssim \exp(t{\g}^{+})$.

From another point of view, we can decompose $\dot u$ and shift the component related to $\Pi_{3}$. Then we define 
\be
U_{1}=\op_{\e}{(\Pi_1)}\dot u, \quad U_{2}=\op_{\e}{(\Pi_2)}\dot u, \quad U_{3}=e^{i\th}\op_{\e}{(\Pi_3)}\dot u, \quad U_{4}=\op_{\e}{(\Pi_4)}\dot u,\quad \cdots, \quad U_{J}=\op_{\e}{(\Pi_J)}\dot u.\nn
\ee
Then the perturbation unknown $\dot u$ is decomposed by
\be
\dot u= U_{1}+U_{2}+e^{-i\th}U_{3}+\cdots+U_{J}\nn
\ee
and
\be
A=\diag(\l_{1},\l_{2},\l_{3,-1}+\omega,\cdots,\l_{J}).\nn
\ee

Similar to Case 1 and Case 2 in Section \ref{Cases of two coupled resonances}, we have
\ba
  B=\sum_{p=\pm 1}
  \bp  e^{ip\th}\Pi_{1,+p}B_{p}\Pi_{1} &  e^{ip\th}\Pi_{1,+p}B_{p}\Pi_{2}  &  e^{i(p-1)\th}\Pi_{1,+(p-1)}B_{p}\Pi_{3,-1}&\cdots& e^{ip\th}\Pi_{1,+p}B_{p}\Pi_{J}
    \\ 
     e^{ip\th}\Pi_{2,+p}B_{p}\Pi_{1} & e^{ip\th}\Pi_{2,+p}B_{p}\Pi_{2} & e^{i(p-1)\th}\Pi_{2,+(p-1)}B_{p}\Pi_{3,-1}&\cdots& e^{ip\th}\Pi_{2,+p}B_{p}\Pi_{J}
    \\
     e^{i(p+1)\th}\Pi_{3,+p}B_{p}\Pi_{1} & e^{i(p+1)\th}\Pi_{3,+p}B_{p}\Pi_{2} &  e^{ip\th}\Pi_{3,+(p-1)}B_{p}\Pi_{3,-1} &\cdots&e^{i(p+1)\th}\Pi_{3,+p}B_{p}\Pi_{J}
 \\
 \vdots &\vdots &\vdots & &\vdots 
 \\
     e^{ip\th}\Pi_{J,+p}B_{p}\Pi_{1} & e^{ip\th}\Pi_{J,+p}B_{p}\Pi_{2} &  e^{i(p-1)\th}\Pi_{J,+(p-1)}B_{p}\Pi_{3,-1} &\cdots&e^{ip\th}\Pi_{J,+p}B_{p}\Pi_{J}\ep.\nn
  \ea
Since we only have $(1,2), (1,3), (2,3)\in \Re_{0}$ which means $(2,1)\notin\Re_{0}$.  Then after normal form reduction,  we have 
\ba\label{(2.40)}
  \check{B}=\bp  0 &  \chi_{12}e^{i\th}b_{12}^{+}  & \chi_{13,-1}(b_{13}^{+})_{-1}
    \\ 
     e^{-i\th}\chi_{12,-1}(b_{21}^{-})_{-1} & 0 & \chi_{23,-1}(b_{23}^{+})_{-1}
     \\
     \chi_{13,-1}(b_{31}^{-})_{-1}  & \chi_{23,-1}(b_{32}^{-})_{-1} &  0 \ep.
  \ea
  
Assume $\mathcal {R}_{12}\bigcap(\mathcal {R}_{13}+k)\bigcap(\mathcal {R}_{23}+k)\neq \varnothing$ and $b_{12}^{+}$ is transparent in $\mathcal {R}_{12}$, or we assume $(\mathcal {R}_{12}+k)\bigcap(\mathcal {R}_{13}+k)\bigcap(\mathcal {R}_{23}+k)\neq\varnothing$ and $b_{21}^{-}$ is transparent in $\mathcal {R}_{12}$, then the oscillation terms $\chi_{12}e^{i\th}b_{12}^{+}$, $e^{-i\th}\chi_{12,-1}(b_{21}^{-})_{-1}$ will be eliminated in normal form reduction. Now we use the assumption  one of  $b_{13}^{+}$, $b_{31}^{-}$ is transparent in $\mathcal {R}_{13}$, or one of  $b_{23}^{+}$, $b_{32}^{-}$ is transparent in $\mathcal {R}_{23}$, the coupled case of $\check{B}$ could be to one of the following cases
\ba
       \check{B}_{1}=\bp  0 &  0  & \chi_{13,-1}(b_{13}^{+})_{-1}
    \\ 
     0 & 0 & \chi_{23,-1}(b_{23}^{+})_{-1}
     \\
     \chi_{13,-1}(b_{31}^{-})_{-1}  & 0 &  0 \ep,
       \check{B}_{2}=\bp  0 &  0  & \chi_{13,-1}(b_{13}^{+})_{-1}
    \\ 
     0 & 0 & 0
     \\
     \chi_{13,-1}(b_{31}^{-})_{-1}  & \chi_{23,-1}(b_{32}^{-})_{-1} &  0 \ep,\nn
  \ea
  \ba
  \check{B}_{3}=\bp  0 &  0  & \chi_{13,-1}(b_{13}^{+})_{-1}
    \\ 
     0 & 0 & \chi_{23,-1}(b_{23}^{+})_{-1}
     \\
    0  & \chi_{23,-1}(b_{32}^{-})_{-1} &  0 \ep,
\check{B}_{4}=\bp  0 &  0  & 0
    \\ 
     0 & 0 & \chi_{23,-1}(b_{23}^{+})_{-1}
     \\
     \chi_{13,-1}(b_{31}^{-})_{-1}  & \chi_{23,-1}(b_{32}^{-})_{-1} &  0 \ep.\nn
  \ea
 Combine Remark \ref{Remark 2.1} and the same method as Case 3  in Section \ref{Cases of two coupled resonances} to obtain the bound of symmetric flows: $|S_{0}(\tau,t)|\lesssim\exp(t{\g}^{+})$. Similar to Section \ref{Instability of WKB solution}, we finally obtain the second result in Theorem \ref{New instability criterion} by the same method of Duhamel representation to establish the instability.

\section{Applications}
\label{Coupled Klein-Gordon systems with equal masses}
In this section, we give examples comprising coupled Klein-Gordon systems in $\mathbb R^{d}$ with equal masses and different velocities to show the application of our theorem. We have the following Klein-Gordon operators 
\be
\partial_{t}+A_{1}(\partial_{x})+\frac{1}{\e}L_{0}, \quad\partial_{t}+A_{1}(\th_{0}\partial_{x})+\frac{1}{\e}L_{0}\nn
\ee
with $0<\th<1$ implying different velocities and
\ba\label{(3.1)}
  A_{1}(\partial_{x})=\bp  0 &  \partial_{x}  & 0
    \\ 
     \partial_{x}\cdot & 0 & 0
     \\
     0 & 0 &  0 \ep,\quad
    L_{0}=\bp  0 &  0  & 0
    \\ 
    0 & 0 & \omega_{0}
     \\
     0 & -\omega_{0} &  0 \ep,
  \ea
where $\omega_{0}>0$ and $x\in \mathbb R^{d}$. For $U=(u,v)=(u_{1},u_{2},u_{3},v_{1},v_{2},v_{3})\in \mathbb R^{2(d+2)}$ with $u_{1}\in \mathbb R^{d}$, $u_{2}\in \mathbb R$, $u_{3}\in \mathbb R$, $v_{1}\in \mathbb R^{d}$, $v_{2}\in \mathbb R$, $v_{3}\in \mathbb R$, we give the coupled systems as follows
\begin{eqnarray}\label{(3.2)}
\begin{cases}
(\partial_{t}+A_{1}(\partial_{x})+L_{0})u=\frac{1}{\sqrt{\e}}B^{1}(U,U),\\
(\partial_{t}+A_{1}(\th_{0}\partial_{x})+L_{0})v=\frac{1}{\sqrt{\e}}B^{2}(U,U),
\end{cases}
\end{eqnarray}
where $B^{1}$ and $B^{2}$ are bilinear $\mathbb R^{2(d+2)}\times \mathbb R^{2(d+2)}\rightarrow \mathbb R^{d+2}$ specified as
\be\label{(3.3)}
 B^{1}(U,U^{'})=\frac{1}{2}\bp 0\\ u_{2}v_{2}^{'}+u_{2}^{'}v_{2}+u_{2}u_{2}^{'}\\0\ep,\quad
B^{2}(U,U^{'})=\frac{1}{2}\bp 0\\ u_{3}v_{3}^{'}+u_{3}^{'}v_{3}+u_{3}u_{3}^{'}\\0\ep
\ee
with $U=(u_{1},u_{2},u_{3},v_{1},v_{2},v_{3})$ and $U^{'}=(u_{1}^{'},u_{2}^{'},u_{3}^{'},v_{1}^{'},v_{2}^{'},v_{3}^{'})$. We denote 
\ba\label{(3.4)}
  A(\partial_{x})=\bp  A_{1}(\partial_{x}) &  0
    \\ 
     0 & A_{1}(\th_{0}\partial_{x})\ep,\quad
    A_{0}=\bp  L_0 & 0
    \\ 
    0 & L_0 \ep,\quad
    B(U,U^{'})=\bp B^{1}(U,U^{'})\\ B^{2}(U,U^{'})\ep.
  \ea
Then, we obtain the 1-multiplicity eigenvalues for matrix $A_{1}(\xi)+L_{0}/i$ 
\ba
\l_{1}(\xi)=\sqrt{\omega_{0}^{2}+|\xi|^{2}}=-\l_{4}(\xi),\quad
\l_{2}(\xi)=\sqrt{\omega_{0}^{2}+\th_{0}^{2}|\xi|^{2}}=-\l_{3}(\xi)\nn
\ea
and 2d-multiplicity eigenvalues $\l_{5}(\xi)=0$. As discussed in \cite{Lu}, corresponding smooth spectral decomposition in Assumption \ref{Smooth spectral decomposition} and \eqref{(1.4)} is satisfied for $A_{1}(\partial_{x})+L_{0}/i$ and $A_{1}(\th_{0}\partial_{x})+L_{0}/i$.
\subsection{Verification of Assumption \ref{WKB solution}: WKB solution}
\label{Verification of Assumption 1.2: WKB solution}
We select the phase $\b=(\o,k)$ belongs to fast Klein-Gordon branch on the variety (i.e. $\o=\sqrt{\o_{0}^{2}+|k|^{2}}$) satisfying 
\begin{itemize}
\item The only harmonic of $\b$ on the fast branches are $p\in \{-1,1\}$, $p^{2}\o^{2}=\o_{0}^{2}+p^{2}|k|^{2}\Rightarrow p\in \{-1,1\}$;
 \item No harmonic of $\b$ belongs to the slow branches i.e. $p^{2}\o^{2}\neq\o_{0}^{2}+p^{2}\th_{0}^{2}|k|^{2}$;
\item No auto-resonances occur i.e. for $\xi\in \mathbb R^{d}$, $\l_{1}(\xi+k)=\pm\o+\l_{1}(\xi)$ and $\l_{2}(\xi+k)=\pm\o+\l_{2}(\xi)$ have no solution;
\end{itemize}
where $\o\in \mathbb R$ is a characteristic temporal frequency and $k\in\mathbb R^{d}$ is the initial wavenumber. 
By \cite{Lu}, it suffices to check the weak transparency condition (6.48) in \cite{Lu} to verify Assumption \ref{WKB solution}. Denoting $\Pi(p\b)$ the orthogonal projector onto $\ker(-ip\o+A(ipk)+A_{0})$, for $|p|=1$, the $\ker(-ip\o+A(ipk)+A_{0})$ are 1-dimensional and generated by $\vec{e}_{1}$ and  $\vec{e}_{-1}=(\vec{e}_{1})^{*}$. Then we obtain 
\be\label{(3.5)}
\vec{e}_{1}=\frac{1}{\sqrt{2}}(\frac{k}{\o},1,\frac{i\o_{0}}{\o},0_{{\mathbb C}^{d+2}})\in {\mathbb C}^{2(d+2)}
\ee
and 
\be\label{(3.6)}
\Pi(\b)U=(U,\vec{e}_{1})\vec{e}_{1},\quad\Pi(-\b)U=(U,\vec{e}_{-1})\vec{e}_{-1},
\ee
where $(\cdot,\cdot)$ represents the Hermitian scalar product in ${\mathbb C}^{2(d+2)}$. For $p=0$, the kernels are generated by
\be
\vec{e}_{0}=\frac{1}{\sqrt{|u_{1}|^{2}+|v_{1}|^{2}}}(u_{1},0,0,v_{1},0,0),\nn
\ee
and the orthogonal projector $\Pi(0)$ onto the kernel of $A(0)+A_{0}/i$ is
\be
\Pi(0)U=(U,\vec{e}_{0})\vec{e}_{0}=(u_{1},0,0,v_{1},0,0).\nn
\ee

Thus by definition of $B$ in \eqref{(3.4)}, for any $U,U^{'}$, we have
\be\label{(3.7)}
\Pi(0)B(U,U^{'})=0, B(\Pi(0)U,U^{'})=0, B(U, \Pi(0)U^{'})=0.
\ee
Then we get
\be
\Pi(0)B(\Pi(0)U,\Pi(0)U^{'})=0,~for~ p=0.\nn
\ee
and for $|p|=1$:
\ba
~~~&\Pi(\b)\sum_{p_{1}+p_{2}=1}B(\Pi(p_{1}\b)U,\Pi(p_{2}\b)U^{'})=\Pi(\b)B(\Pi(\b)U,\Pi(0\cdot\b)U^{'})+\Pi(\b)B(\Pi(0\cdot\b)U,\Pi(\b)U^{'})=0,\\
&\Pi(-\b)\sum_{p_{1}+p_{2}=-1}B(\Pi(p_{1}\b)U,\Pi(p_{2}\b)U^{'})=\Pi(-\b)B(\Pi(-\b)U,\Pi(0\cdot\b)U^{'})+\Pi(-\b)B(\Pi(0\cdot\b)U,\Pi(-\b)U^{'})=0.\nn
\ea

Then we verified (6.48) in \cite{Lu} and furthermore, we obtain the existence of WKB approximate solution.
\subsection{Verification of Assumption \ref{Relaxed  separation conditions}: Resonances and transparency}
\label{Verification of Relaxed  separation conditions}
Combining the form of the characteristic variety and the $\b$ chosen in Section \ref{Verification of Assumption 1.2: WKB solution}, the resonant pairs are
\be\label{(3.8)}
\Re=\{(1,2),(1,5),(2,5),(3,4),(5,3),(5,4)\}.
\ee
Then we have
\be\label{(3.9)}
B(\vec{e}_{1})U=B(\vec{e}_{1},U)+B(U,\vec{e}_{1})=\frac{1}{\sqrt{2}}(0_{d},u_{2}+v_{2},0,0_{d},\frac{i\o}{\o}(u_{3}+v_{3}),0),
\ee
\be\label{(3.10)}
B(\vec{e}_{-1})U=B(\vec{e}_{-1},U)+B(U,\vec{e}_{-1})=\frac{1}{\sqrt{2}}(0_{d},u_{2}+v_{2},0,0_{d},-\frac{i\o}{\o}(u_{3}+v_{3}),0).
\ee
We denote $U_{5}(\xi)$ a element in the image of $\Pi_{5}(\xi)$ and $U_{5}(\xi)$ satisfies the following form
\be\label{(3.11)}
U_{5}(\xi)=(u_{1},0,\frac{-i\xi\cdot u_{1}}{\o_{0}},v_{1},0,\frac{-i\th_{0}\xi\cdot u_{1}}{\o_{0}}),\quad  (u_{1},v_{1})\in \mathbb C^{d}\times\mathbb C^{d}.
\ee
After a brief observation, we have
\be\label{(3.12)}
\Pi_{5}(\cdot)B(\vec{e}_{\pm1})=0,
\ee
which means for any $U\in \mathbb C^{2(d+2)}$, $B(\vec{e}_{\pm1})U$ belongs to the orthogonal of the range of $\Pi_{5}(\xi)$.

The other projectors are
\be
\Pi_{j}(\xi)U=(U,\O_{j}(\xi))\O_{j}(\xi), \quad j=1,2,3,4,\nn
\ee
where
\be
\O_{j}(\xi)=\frac{1}{\sqrt{2}}(\frac{\xi}{\l_{j}},1,\frac{i\o_{0}}{\l_{j}},0_{\mathbb C^{d+2}}), \quad j=1,4;\nn
\ee
\be
\O_{j^{'}}(\xi)=\frac{1}{\sqrt{2}}(0_{\mathbb C^{d+2}},\frac{\th_{0}\xi}{\l_{j^{'}}},1,\frac{i\o_{0}}{\l_{j^{'}}}), \quad j^{'}=2,3.\nn
\ee
Then we have
\be\label{(3.13)}
(\O_{1}(\xi+k),B(\vec{e}_{1})U_{5}(\xi))=0,
\ee
\be\label{(3.14)}
(\O_{4}(\xi),B(\vec{e}_{-1})U_{5}(\xi+k))=0.
\ee

Combining $\Pi_{5}(\cdot)B(\vec{e}_{\pm1})=0$, it shows that the resources $(1,5),(5,4)$ are transparent and furthermore, we have
\be
(\O_{1}(\xi+k),B(\vec{e}_{1})\O_{2}(\xi))=\frac{1}{2\sqrt{2}},\nn
\ee
\be
(\O_{2}(\xi),B(\vec{e}_{-1})\O_{1}(\xi+k))=\frac{1}{2\sqrt{2}}\cdot\frac{\o_{0}^{2}}{\o\l_{1}(\xi+k)},\nn
\ee
\be
(\O_{3}(\xi+k),B(\vec{e}_{1})\O_{4}(\xi))=-\frac{1}{2\sqrt{2}}\cdot\frac{\o_{0}^{2}}{\o\l_{4}(\xi)},\nn
\ee
\be
(\O_{4}(\xi),B(\vec{e}_{-1})\O_{3}(\xi+k))=\frac{1}{2\sqrt{2}},\nn
\ee
\be
(\O_{2}(\xi+k),B(\vec{e}_{1})U_{5}(\xi))=\frac{\xi\cdot u_{1}+\th_{0}\xi\cdot v_{1}}{2\o},\nn
\ee
\be
(\O_{3}(\xi),B(\vec{e}_{-1})\O_{5}(\xi+k))=-\frac{(\xi+k)(u_{1}+\th_{0}v_{1})}{2\o}.\nn
\ee
Thus we get
\be\label{(3.15)}
\Re_{0}=\{(1,2),(2,5),(3,4),(5,3)\},
\ee
in which set no auto-resonances are contained. Furthermore, we consider the following corresponding resonant sets:
\be
\mathcal R_{12}=\{\xi: ~\sqrt{\o_{0}^{2}+|\xi+k|^{2}}=\o+\sqrt{\o_{0}^{2}+\th_{0}^{2}|\xi|^{2}}\},\nn
\ee
\be
\mathcal R_{25}=\{\xi: ~\th_{0}|\xi+k|=|k|\},\nn
\ee
\be
\mathcal R_{34}=\{\xi: ~\sqrt{\o_{0}^{2}+|\xi|^{2}}-\sqrt{\o_{0}^{2}+\th_{0}^{2}|\xi+k|^{2}}=\o\},\nn
\ee
\be
\mathcal R_{53}=\{\xi: ~\th_{0}|\xi|=|k|\},\nn
\ee
where $\o=\sqrt{\o_{0}^{2}+|k|^{2}}$. Now, we are trying to verify separation conditions in Assumption \ref{Relaxed  separation conditions}. 

Firstly, we consider $(1,2)$-non-transparent resonance and its coupled $(2,5)$-resonance in $\Re_{0}$ and give the following frequency shift:
\be
U_{1}=e^{i\th}\op_{\e}{(\Pi_1)}\dot u, \quad U_{2}=\op_{\e}{(\Pi_2)}\dot u, \quad U_{3}=\op_{\e}{(\Pi_3)}\dot u,\quad U_{4}=\op_{\e}{(\Pi_4)}\dot u,\quad U_{5}=e^{-i\th}\op_{\e}{(\Pi_5)}\dot u.\nn
\ee
Obviously \eqref{(3.12)} shows that $b_{52}^{-}$ is transparent in $\mathcal R_{25}$ i.e. $(b_{52}^{-})_{-1}$ is transparent in $(\mathcal R_{25}+k)$. After normal form reduction, we have 
  \ba\label{(3.16)}
  \check{B}=
  \bp  0 & \chi_{12}b_{12}^{+}  & 0
    \\ 
     \chi_{12}b_{21}^{-} & 0 & \chi_{25,-1}(b_{25}^{+})_{-1}
    \\
    0   & 0 &  0 \ep, 
  \ea
  and
\be
\mathcal R_{12}\cap(\mathcal R_{25}+k)=\{\xi:~|\xi+k|=\sqrt{3\o_{0}^{2}+4|k|^{2}},~\th_{0}|\xi|=|k|\}.\nn
\ee
Since $\o_{0}$ and $|k|$ are fixed, bounded constants and by inequality 
\be
(\frac{1}{\th_{0}}-1)|k|=||\xi|-|k||\lesssim|\xi+k|=\sqrt{3\o_{0}^{2}+4|k|^{2}}=C(\o,k)\lesssim||\xi|+|k||=(\frac{1}{\th_{0}}+1)|k|,\nn
\ee
we can easily draw a conclusion that when $\th_{0}$ small enough, it gives
\be
\mathcal R_{12}\cap(\mathcal R_{25}+k)\neq\varnothing.\nn
\ee
Then in $\mathcal R_{12}\cap(\mathcal R_{25}+k)$, \eqref{(3.16)} can be reduced to
  \ba\label{(3.17)}
  \check{B}=\chi_{12}\chi_{25,-1}
  \bp  0 & b_{12}^{+}  & 0
    \\ 
     b_{21}^{-} & 0 & (b_{25}^{+})_{-1}
    \\
    0   &0 &  0 \ep.
  \ea
 By Case 1, we get 
$|S_{0}(\tau,t)|\lesssim\exp(t\g_{12}^{+})$ for corresponding symbolic flow. When $\xi\in\mathcal R_{12}\backslash(\mathcal R_{25}+k)$ and $\xi\in(\mathcal R_{25}+k)\backslash\mathcal R_{12}$, we will obtain 
$|S_{0}(\tau,t)|\lesssim \exp(t\g_{12}^{+})$ and $|S_{0}(\tau,t)|\lesssim 1$ respectively. When $\th_{0}$ does not satisfy the smallness and $\mathcal R_{12}\cap(\mathcal R_{25}+k)=\varnothing$, we can obtain the same bounds for symbolic flow both $\xi\in\mathcal R_{12}$ and $\xi\in(\mathcal R_{25}+k)$ .

Then, we consider $(2,5)$-non-transparent resonance and its coupled $(1,2)$-resonance and $(5,3)$-resonance in $\Re_{0}$. We make the frequency shift as 
\be
U_{1}=e^{-2i\th}\op_{\e}{(\Pi_1)}\dot u,\quad U_{2}=e^{-i\th}\op_{\e}{(\Pi_2)}\dot u,\quad U_{3}=e^{i\th}\op_{\e}{(\Pi_3)}\dot u,\quad U_{4}=\op_{\e}{(\Pi_4)}\dot u,\quad U_{5}=\op_{\e}{(\Pi_5)}\dot u,\nn
\ee
where $\th=(k\cdot x-wt)/{\e}$. Combining (3.12), we know that $b_{52}^{-}$ and $(b_{53}^{+})_{-1}$ are transparent in $\mathcal R_{25}$ and $(\mathcal R_{53}+k)$ respectively, then we have
{\ba\label{(3.18)}
  \check{B}=
  \bp  0 & \chi_{12,+1}(b_{12}^{+})_{+1}  &  0&0
    \\ 
    \chi_{12,+1}(b_{21}^{-})_{+1} & 0 &0&\chi_{25}b_{25}^{+}
   \\
    0  & 0 &  0&\chi_{53,-1}(b_{35}^{-})_{-1}
   \\
    0   & 0 & 0&0
    \ep.
  \ea}
Since
\be
\mathcal R_{25}\cap(\mathcal R_{53}+k)=\{\xi:~|\xi|=\sqrt{1-\th_{0}^{2}}|k|\},\nn
\ee
\be
\mathcal R_{25}\cap(\mathcal R_{12}-k)=\{\xi:~|\xi+2k|=\sqrt{3\o_{0}^{2}+4|k|^{2}},~\th_{0}|\xi+k|=|k|\},\nn
\ee
they show that 
\be
\mathcal R_{25}\cap(\mathcal R_{53}+k)\cap(\mathcal R_{12}-k)=\{\xi:~|\xi|=\sqrt{1-\th_{0}^{2}}|k|,~|\xi|^{2}=3\o_{0}^{2}\}.\nn
\ee
If the fixed $\xi$ and $\o_{0}$ satisfy $(1-\th_{0}^{2})|k|=3\o_{0}^{2}$, then $\mathcal R_{25}\cap(\mathcal R_{53}+k)\cap(\mathcal R_{12}-k)\neq\varnothing$, in which intersection we have 
 {\ba\label{(3.19)}
  \check{B}=\chi_{25}\chi_{53,-1}\chi_{12,+1}
  \bp  0 & (b_{12}^{+})_{+1}  &  0&0
    \\ 
    (b_{21}^{-})_{+1} & 0 &0&b_{25}^{+}
   \\
    0  & 0 &  0&(b_{35}^{-})_{-1}
   \\
    0   &0 &  0&0
    \ep.
  \ea}
 Then using similar estimates in Case 1 in Section \ref{Cases of two coupled resonances}, we have $|S_{0}(\tau,t)|\lesssim \exp(t\g_{12}^{+})$.

For the fixed $\xi$ and $\o_{0}$, if  $(1-\th_{0}^{2})|k|\neq3\o_{0}^{2}$, then $\mathcal R_{25}\cap(\mathcal R_{53}+k)\cap(\mathcal R_{12}-k)=\varnothing$, then we only need to discuss the case of intersection $\mathcal R_{25}\cap(\mathcal R_{53}+k)$ and $\mathcal R_{25}\cap(\mathcal R_{12}-k)$, in those intersections we have
  \ba\label{(3.20)}
  \check{B_{1}}=\chi_{25}\chi_{53,-1}
  \bp  0 & 0  &  b_{25}^{+}
    \\ 
     0 & 0 &(b_{35}^{-})_{-1}
    \\
    0   &0 &  0 \ep,\quad
     \check{B_{2}}=\chi_{25}\chi_{12,+1}
  \bp  0 & (b_{12}^{+})_{+1}  &  0
    \\ 
    (b_{21}^{-})_{+1} & 0&b_{25}^{+}
   \\
    0   & 0 &0
    \ep.
  \ea
Then we can give the bounds $|S_{0}(\tau,t)|\lesssim 1$ and $|S_{0}(\tau,t)|\lesssim \exp(t\g_{12}^{+})$ respectively for  symbolic flow directly by same procedures in Case 1 of Section \ref{Cases of two coupled resonances}.

Furthermore,  we consider $(3,4)$-non-transparent resonance and its coupled $(5,3)$-resonance in $\Re_{0}$. We make the frequency shift as
\be
U_{1}=\op_{\e}{(\Pi_1)}\dot u, \quad U_{2}=\op_{\e}{(\Pi_2)}\dot u, \quad U_{3}=\op_{\e}{(\Pi_3)}\dot u, \quad U_{4}=e^{i\th}\op_{\e}{(\Pi_4)}\dot u, \quad U_{5}=e^{-i\th}\op_{\e}{(\Pi_5)}\dot u,\nn
\ee
and 
\be
(\mathcal R_{34}+k)\cap\mathcal R_{53}=\{\xi:~|\xi-k|=\sqrt{3\o_{0}^{2}+4|k|^{2}},~\th_{0}|\xi|=|k|\},\nn
\ee
then combining (3.12), we know that $b_{53}^{+}$ are transparent in $\mathcal R_{53}$, we have
  \ba
  \check{B}=\chi_{34,-1}\chi_{53}
  \bp  0 & (b_{34}^{+})_{-1}  & b_{35}^{-}
    \\ 
    ( b_{43}^{-})_{-1} & 0 &0
    \\
   0  & 0 &  0 \ep. \nn
  \ea
Then we can directly use the theory in Case 1 in Section \ref{Cases of two coupled resonances} to give the discussion for the bound of symbolic flow: $|S_{0}(\tau,t)|\lesssim \exp(t\g_{34}^{+})$.

Finally, we consider $(5,3)$-non-transparent resonance and its coupled $(2,5)$-resonance and  $(3,4)$-resonance in $\Re_{0}$. We make the frequency shift as 
\be
U_{1}=\op_{\e}{(\Pi_1)}\dot u, \quad U_{2}=e^{-2i\th}\op_{\e}{(\Pi_2)}\dot u, \quad U_{3}=\op_{\e}{(\Pi_3)}\dot u, \quad U_{4}=e^{i\th}\op_{\e}{(\Pi_4)}\dot u, \quad U_{5}=e^{-i\th}\op_{\e}{(\Pi_5)}\dot u,\nn
\ee
We know that $(b_{52}^{-})_{+1}$ and  $b_{53}^{+}$ are transparent in $(\mathcal R_{25}-k)$ and $\mathcal R_{53}$ respectively, thus we have
  \ba\label{(3.21)}
  \check{B}=
  \bp  0 & 0  &  0&\chi_{25,+1}(b_{25}^{+})_{+1}
   \\ 
   0 & 0 &\chi_{34,-1}(b_{34}^{+})_{-1}&\chi_{53}b_{35}^{-}
    \\
   0  &  \chi_{34,-1}(b_{43}^{-})_{-1} &  0&0
    \\
   0  & 0&0&0
     \ep.
  \ea

We need to consider the following intersections 
\be
\mathcal R_{53}\cap(\mathcal R_{25}-k)=\{\xi:~|\xi+k|=\sqrt{1-\th_{0}^{2}}|k|\},\nn
\ee
\be
\mathcal R_{53}\cap(\mathcal R_{34}+k)=\{\xi:~|\xi-k|=\sqrt{3\o_{0}^{2}+4|k|^{2}},~\th_{0}|\xi|=|k|\}.\nn
\ee
Moreover, we have 
\be
\mathcal R_{53}\cap(\mathcal R_{25}-k)\cap(\mathcal R_{34}+k)=\{\xi:~|\xi|=\frac{1}{\th_{0}}|k|,~|\xi|^{2}=3\o_{0}^{2}+|k|^{2}\},\nn
\ee

If the fixed $\xi$, $\o_{0}$ and $\th_{0}$ satisfy $\frac{1}{\th_{0}^{2}}|k|^{2}=3\o_{0}^{2}+|k|^{2}$, then $\mathcal R_{53}\cap(\mathcal R_{25}-k)\cap(\mathcal R_{34}+k)\neq\varnothing$, in which intersection we have 
  \ba\label{(3.22)}
  \check{B}=\chi_{53}\chi_{25,+1}\chi_{34,-1}
  \bp  0 & 0  &  0&(b_{25}^{+})_{+1}
   \\ 
   0 & 0 &(b_{34}^{+})_{-1}&b_{35}^{-}
    \\
   0  &  (b_{43}^{-})_{-1} &  0&0
    \\
   0   & 0&0&0
     \ep.
  \ea
then we have $|S_{0}(\tau,t)|\lesssim \exp(t\g_{34}^{+})$.
If we assume  $\frac{1}{\th_{0}^{2}}|k|^{2}\neq3\o_{0}^{2}+|k|^{2}$, then $\mathcal R_{53}\cap(\mathcal R_{25}-k)\cap(\mathcal R_{34}+k)=\varnothing$, in both intersections $\mathcal R_{53}\cap(\mathcal R_{25}-k)$ and $\mathcal R_{53}\cap(\mathcal R_{34}+k)$, \eqref{(3.21)} is converted to
\ba\label{(3.23)}
  \check{B_{1}}=\chi_{53}\chi_{25,+1}
  \bp  0 & 0  & (b_{25}^{+})_{+1}
   \\ 
   0 & 0 &b_{35}^{-}
    \\
   0   & 0&0
     \ep,\quad
     \check{B_{2}}=\chi_{53}\chi_{34,-1}
     \bp 
    0 &(b_{34}^{+})_{-1}&b_{35}^{-}
    \\
     (b_{43}^{-})_{-1} &  0&0
    \\
   0&0&0
     \ep.
  \ea
Thus the upper bounds of symbolic flow are obtained as: $|S_{0}(\tau,t)|\lesssim\exp(t\g^{+})$ and we can apply Theorem \ref{New instability criterion} to coupled Klein-Gordon systems with equal masses whose nonlinear terms specified as \eqref{(3.3)} for $x,\xi\in \mathbb R^{d}(d\geq 1)$. After a short observation, we know: $\Gamma_{12}(\xi)>0\Rightarrow {\Gamma}>0$ and combining Theorem \ref{New instability criterion} (the relaxed instability criterion), we obtain the instability of WKB solutions $u_{a}$. 

\begin{remark}\label{Remark 3.1} In Section 5.2 of \cite{Lu}, under the Assumption 2.8 (Assumption \ref{Separation conditions} in this paper), the example 'coupled Klein-Gordon systems with equal masses' with different nonlinear term $B(U,V)$ defined in \cite{Lu} in (5.39) shows the instability criterion holds only in 1-dimensional space. Now we are trying to verify the relaxed instability criterion in Theorem \ref{New instability criterion} to improve the results i.e. to turn the restriction 1-dimensional space into d dimensional space with any $d\geq 1$. 

For $\Re_{0}=\{(1,2),(1,5),(3,4),(5,4)\}$, we already have 
\begin{itemize}
\item$\mathcal R_{12}\cap\mathcal R_{15}=\varnothing$ for $(1,2)$-resonance and its coupled $(1,5)$-resonance in $\Re_{0}$;
\item$\mathcal R_{15}\cap\mathcal R_{12}\cap(\mathcal R_{54}+k)=\varnothing$ for $(1,5)$-resonance and its coupled $(1,2)$, $(5,4)$-resonance in $\Re_{0}$;
\item$\mathcal R_{34}\cap\mathcal R_{54}=\varnothing$ for $(3,4)$-resonance and its coupled $(5,4)$-resonance in $\Re_{0}$;
\item$\mathcal R_{54}\cap(\mathcal R_{15}-k)\cap\mathcal R_{34}=\varnothing$ for $(5,4)$-resonance and its coupled $(1,5)$, $(3,4)$-resonance in $\Re_{0}$.
\end{itemize}
Combining (5.41) in \cite{Lu}: $\Pi_{5}(\cdot)B(\vec{e}_{\pm 1})=0$, we know that $b_{51}^{-}$  and $ (b_{54}^{+})_{-1}$ are transparent in $\mathcal R_{15}$ and $(\mathcal R_{54}+k)$ respectively. Thus we only need to consider $\mathcal R_{15}\cap(\mathcal R_{54}+k)=\{\xi\in\mathbb R^{d}:~|\xi|=0\}$, in which intersection by the method of frequency shift and normal form reduction, we have 
  \ba\label{(3.24)}
  \check{B}=\chi_{15}\chi_{54,-1}
  \bp  0 & 0  & b_{15}^{+}
    \\ 
    0 & 0 &(b_{45}^{-})_{-1}
    \\
    0  & 0 &  0 \ep, 
  \ea
Then we obtain $|S_{0}(\tau,t)|\lesssim 1$ directly by Case 1 in Section \ref{Cases of two coupled resonances}. Hence Theorem \ref{New instability criterion} (the relaxed instability criterion) is applicable for this case. Combining the calculation result: $\Gamma_{12}(\xi)>0\Rightarrow {\Gamma}>0$, it shows instability for coupled Klein-Gordon systems with equal masses in d dimensional space whose nonlinear terms specified as (5.39) in \cite{Lu}. 
\end{remark}

{\bf Acknowledgment}: J. Pan is partially supported by the NSF of China under Grant 12171235.


\begin{thebibliography}{000}
\bibitem{Cheverry} C. Cheverry, O. Gue\`s, G. M\'etivier. Oscillations fortes sur un champ lin\'eairement d\'eg\'en\'er\'e, \emph{Ann. Sci. E\'cole Normale Sup.} 36 (2003) 691-745.
\bibitem{Dumas} E. Dumas. About nonlinear geometric optics, \emph{Bol. Soc. Esp. Mat. Apl. SeMA.} 35 (2006) 7-41.
\bibitem{Dumas1} E. Dumas. Diffractive optics with curved phases: beam dispersion and transitions between light and shadow, \emph{Asymptotic Anal.} 38 (2004) 47-91.
\bibitem{Joly} J.--L. Joly, G M\'etiver, J. Rauch. Transparent nonlinear geometric optics and Maxwell-Bloch equations, \emph{ J. Differential Equations.} 166 (2000) 175-250.
\bibitem{Joly1}  J.--L. Joly, G M\'etiver, J. Rauch. Coherent and focusing multidimensional nonlinear geometric optics, \emph{Ann. Sci. E\'cole Norm. Sup.} 28 (1995) 51-113.
\bibitem{Jeanne} P.--Y. Jeanne. Optique g\'eom'etrique pour des syst\`smes semi-lin\'eaires avec invariance de jauge, \emph{M\'emoire Soc. Math. Fr.} 90 (2002). 
\bibitem{Klainerman} S. Klainerman. The null condition and global existence to nonlinear wave equations, \emph{ Lectures in Appl. Math.} 23 (1986).   \bibitem{Lannes} D.  Lannes. Space time resonances, after Germain, Masmoudi, Shatah. \emph{Asterisque-Societe Mathematique de France.} 352 (2013).
\bibitem{Lu} Y. Lu, B Texier. A Stability Criterion for High-Frequency Oscillations, \emph{M\'emoires de la Soci\'et\'e Math\'ematique de France.} 1 (2013).
\bibitem{Lu1} Y. Lu. High-frequency limit of the Maxwell-Landau-Lifshitz system in the diffractive optics regime, \emph{Asymptotic Anal.} 82 (2013) 109-137.
\bibitem{Lu2}Y. Lu, Z. Zhang. Partially strong transparency conditions and a singular localization method in geometric optics, \emph{Arch. Ration. Mech. Anal.} 222 (1) (2016) 245-283.
\bibitem{Lu3}Y. Lu. Higher-order resonances and instability of high-frequency WKB solutions. \emph{J. Differ. Equations.} 260 (3) (2016) 2296-2353.
\bibitem{Lerner} N. Lerner, Y. Morimoto, C.--J. Xu, Instability of the Cauchy-Kovalevskaya solution for a class of nonlinear systems, \emph{Amer. J. Math.} 132 (2010) 99-123.
\bibitem{Lerner1} N. Lerner, T. Nguyen, B. Texier. The onset of instability in first-order systems. \emph{J. Eur. Math. Soc.} 20 (6) (2015).
\bibitem{Metivier} G. M\'etivier. Para-differential Calculus and Applications to the Cauchy Problem for Nonlinear Systems, Centro di Ricerca Matematica Ennio De Giorgi, CRM Series, vol. 5, Edizioni della Normale, Pisa, 2008.
\bibitem{Nishiura}  Y. Nishiura, M. Mimura. Layer oscillations in reaction-diffusion systems, \emph{SIAM J. Appl. Math.} 49 (1989) 481-514.
\bibitem{Texier1} B. Texier. Approximations of pseudo-differential flows. \emph{Indiana Univ. Math. J.} 65 (1). (2014).
\bibitem{Texier2} B. Texier. The short wave limit for nonlinear, symmetric hyperbolic systems, \emph{Adv. Differ. Equ.} 9 (2004) 1-52.
\bibitem{Texier4} B. Texier. Derivation of the Zakharov equations, \emph{Arch. Ration. Mech. Anal. }184 (2007) 121-183.
\bibitem{Texier3} B. Texier. WKB asymptotics for the Euler-Maxwell equations, \emph{Asymptot. Anal. }42 (2005) 211-250.
\end{thebibliography}
\end{document}